\documentclass[red,11pt,a4paper]{article}
\usepackage{cite}
\usepackage{amsmath}
\usepackage{amscd}
\usepackage{amssymb}
\usepackage[pdftex]{color,graphicx}
\usepackage[pdfencoding=auto, psdextra]{hyperref}
\usepackage{latexsym}
\usepackage{color}
\usepackage{cases}
\usepackage{xcolor}
\usepackage{comment}
\usepackage[normalem]{ulem}
\usepackage{enumitem}
\usepackage{geometry}

\usepackage[T1]{fontenc}
\usepackage[utf8]{inputenc}

\bf
\usepackage{ika}
\usepackage{mathrsfs}

\newcommand{\dv}{\Div}
\newcommand{\sgn}{{\rm sgn}}
\newcommand{\iintO}[1]{\int_{\mathbb{T}^{6}} #1 \ \dx\, \dy}
\newcommand{\T}{\mathbb{T}}
\newcommand{\Td}{\mathbb{T}^3}

\newcommand{\iintTO}[1]{\int_0^T\!\!\!\! \int_{\mathbb{T}^{6}} #1 \ \dx\, \dy\, \dt}

\newcommand{\iintOM}[1]{ \int_{\T^3} #1 \ \dx\,}
\newcommand{\iintTOM}[1]{\int_0^T\!\!\!\! \int_{\T^3} #1 \ \dx\, \dt}

\title{Construction of weak solutions to the equations of a compressible viscous model}

\author{Nilasis Chaudhuri$^1$, Piotr B. Mucha$^1$, Milan Pokorn\'y$^{2}$}
\date{\today}

\begin{document}
	\maketitle
	\medskip
\centerline{$^1$ Faculty of Mathematics, Informatics and Mechanics,  University of Warsaw,}
	\centerline{Institute of Applied Mathematics and Mechanics, ul. Banacha 2, 02-097 Warszawa, Poland}
	\vspace{5mm}
\centerline{$^2$  Faculty of Mathematics and Physics, Charles University -- Prague,}
\centerline{Mathematical Institute of Charles University,
Sokolovsk\'a 83, 186 75 Praha, Czech Republic.}

\begin{abstract}
    The paper aims on the construction of weak solutions to equations of a model of compressible viscous fluids, being a simplification of the classical compressible Navier-Stokes system. 
    We present a novel scheme for approximating systems that preserves structural integrity by avoiding classical regularization with ${\bf{ -\epsilon \Delta \vr}}$, thus maintaining the transport character of  the continuity equation. Our approach, which necessitates specific conditions on the constitutive equation, accommodates physically relevant models such as isentropic and van der Waals gases, and globally handles non-monotone pressures. From an analytical perspective, our method synthesizes techniques from Feireisl–Lions and Bresch–Jabin to demonstrate the convergence of approximate densities using compensated compactness techniques. We also apply renormalization of the continuity equations and utilize weight techniques to manage unfavorable terms.
\end{abstract}
\medbreak
\medbreak

	\tableofcontents

	\section{Introduction }

The construction of weak solutions to systems of partial differential equations describing the flow of compressible viscous fluids is a well-established topic in the classical theory of mathematical fluid mechanics. Over the last three decades, significant progress has been made, particularly through the pioneering ideas of P. L. Lions \cite{Lions2} and E. Feireisl \cite{F}. More recently, D. Bresch and P.-E. Jabin \cite{BJ} have introduced innovative compactness techniques that have rejuvenated the theory.

However, there is a consensus that the current state of the art for weak solutions to the compressible Navier–Stokes system remains unsatisfactory. The core challenge is to identify an approximate system for the original equations in which the existence of suitably regular solutions is guaranteed. These approximate solutions must adhere to energy-type estimates that align with standard kinetic and potential energy. The main difficulty lies in passing to the limit with these approximate solutions. Given that we only have uniform bounds in $L^p$ spaces for the densities, we achieve merely weak convergence. Due to the presence of nonlinear terms, this weak convergence is insufficient to characterize the limit. Therefore, by utilizing the toolbox of compensated compactness techniques, the objective is to obtain strong or pointwise convergence rather than just weak convergence. While simple in its statement, this aim poses a significant challenge.

A common theme in all the aforementioned seminal results is the approximation of the compressible Navier–Stokes system through the parabolic/viscous regularization of the continuity equation:
\begin{equation}\label{cns-visc}
    \begin{array}{l}
         \vr_t +\Div(\vr \vc{u})
         {\red{\bf -\epsilon \Delta \vr}} =0  \\
         (\vr \vc {u})_t + \Div (\vr \vc{u} \otimes \vc{u}) - \Div \mathbf{S} + \nabla \pi(\vr)=0. 
    \end{array}
\end{equation}
In \eqref{cns-visc}, $\vr$ and $\vc{u}$ stand for the density and velocity of the fluid, respectively; $\pi(\vr)$ represents the pressure, and $\mathbf{S}$ is the viscous stress tensor, which typically depends linearly on $\nabla \vc{u}$. 

This strong regularization of the continuity equation introduces pronounced parabolicity to the approximation but compromises the transport character of the first equation. To be precise, the system requires further modification to feasibly achieve the stated goals. As $\ep$ approaches zero, we recover the original compressible Navier–Stokes system. However, establishing this rigorously within the mathematical framework is a challenging and intricate process. Additionally, to obtain a suitable energy inequality, we must modify the momentum equation.



The presence of ${\red{\bf -\epsilon \Delta \vr}}$ introduces additional challenges, particularly regarding the required integrability of the pressure. However, techniques developed over the past years have successfully addressed these issues for the classical system (see \cite{FNP}, \cite{F}, \cite{NS}, \cite{FN}, \cite{NP1}). For more complex systems, such as multi-fluid models, constructing such approximations is considerably more challenging, as noted by Bresch, Mucha, and Zatorska \cite{BMZ} and Novotn\'y and Pokorn\'y\cite{NP}. Note that another approach, based on construction of approximate solutions via a numerical method was addressed in \cite{FKP}.

This motivates the goal of our manuscript: to find an alternative method for constructing systems of the form \eqref{cns-visc}. Since our focus is on methods related to the convergence of the density, we simplify the Navier–Stokes system by considering a version of the active scalar system, described as follows:
\begin{align}
		&\pt \vr + \Div(\vr \vu)=0,\label{bPF:1}\\
		&\pt \vu +\nabla \pi(\vr) = \Delta \vu,\label{bPF:2}
	\end{align}
	where $\vr $ is the density, $\vu$ is the velocity and $\pi$ is the pressure. The system can be considered in the time interval $(0,T)$ and in the spatial variable in the whole space $\R^d$, in the torus $\T^d$ or in a bounded domain with some boundary conditions. To simplify the presentation, we stay in the torus and fix $d=3$.

We introduce a novel scheme for the approximate system that preserves the system's structure. Specifically, we avoid the classical regularization by ${\red{\bf -\epsilon \Delta \vr}}$, maintaining the transport equation for the density. This approach requires specific conditions on the constitutive equation while preserving physically relevant models, such as isentropic and van der Waals gases. In practice, this method accommodates non-monotone pressures.

From an analytical perspective, our approach synthesizes the techniques of Feireisl–Lions and Bresch–Jabin, demonstrating the convergence of approximate densities using the toolbox of compensated compactness. Our novel technique relies on the renormalization of the continuity equations, with the final step employing weight techniques to eliminate terms with unfavorable signs.

We aim to make our paper self-contained, providing comprehensive proofs and frequently referencing classical methods, such as those based on the oscillation defect measure or the Kolmogorov convergence criterion.

\section{Assumptions, Main results and Idea of the proof}

Recall we consider the simplified problem \eqref{bPF:1}--\eqref{bPF:2} in $(0,T)\times \T^3$ equipped with the corresponding initial conditions for the functions $\vu$ and $\vr$. Moreover, we need to specify the pressure function $\pi(\vr)$.
The goal is to consider the most possible general pressure law.
   We take a smooth non-negative function $\pi\in C^1([0,\infty))\cap C^2((0,\infty))$ such that $\pi(\vr) >0$ for $\vr >0$ and 
\begin{equation}\label{pre1}
	\pi(0)=0 \text{  \ \  and  \ \ }  a_2 \vr^\gamma -C \leq \pi(\vr) \leq C + a_1 \vr^\gamma \mbox{ \ \ with \ \ } \gamma >1,\; a_1,a_2,C>0.
\end{equation}
Additionally, we require 
\begin{equation}\label{pre2}
	|\pi'(\vr)|\leq C \vr^{\gamma-1}, \qquad  |\pi''(\vr)|\leq C \vr^{\gamma-2} \mbox{ \ \ for \ \ } \vr > 1.
\end{equation}
In this case, we define the pressure potential as $\displaystyle  \Pi(\vr)= \vr \int_{1}^{\vr} \frac{\pi(\xi)}{\xi^2} \text{d}\xi$.

Within our analysis we  use more specified forms of $\pi$. Typically, we consider $\pi(\vr) \sim p(\vr) + perturbation$, where $p(\vr)$ is expected to have fine properties like the most classical case $\vr^\gamma$. 

 \subsection{Definition of weak solution}
The problem is supplemented with initial data $(\vr_0, \vu_0)$  such that $\vr_0\geq 0$ a.e. in $\T^3$ and 
\begin{align}\label{id}
E_0(\vr,\vu) =	\int_{\T^3} \left(\frac{1}{2}|\vu_0|^2+ \Pi(\vr_0) \right) \dx <\infty .
\end{align} 
At first, we give the definition of a weak solution.
\begin{df} \label{d:wsol}
	Let the initial data $(\vr_0,\vu_0)$ satisfy \eqref{id} and let the pressure follow \eqref{pre1}--\eqref{pre2}. Then we say $(\vr,\vu)$ is a weak solution to the system \eqref{bPF:1}--\eqref{bPF:2} in the class
	\begin{align*}
		&0\leq \vr \in  C_{\text{weak}}([0,T]; L^{\gamma}(\T^3)) \text{ and } \vu \in  L^2(0,T; W^{1,2}(\T^3) )  \cap C([0,T]; L^{2}(\T^3) ) 
	\end{align*}
	if the following holds: 
	\begin{itemize}[leftmargin=*]
		\item Continuity equation holds
		\begin{align*}
			&\intTO{ \Big(\vr \partial_t \psi + \vr \vu \cdot \nabla \psi\Big)} = -\intO{\vr_0 \psi(0)} ,
		\end{align*}
		for any $\psi \in  C_c^1([0,T)\times \T^d)$. 
		\item Renormalized continuity equation holds \begin{align*}
			&\intTO{ \Big(b(\vr) \partial_t \psi + b(\vr) \vu \cdot \nabla \psi + (b(\vr)-b'(\vr)\vr) \Div \vu \psi \Big)} = -\intO{b(\vr_0) \psi(0)} ,
		\end{align*}
		for any $\psi \in  C_c^1([0,T)\times \T^d) $ and $b\in C^1([0,\infty))$ such that $b'(z) = 0$ for  $z>M_0$ for some $M_0 >0$.
		\item The equation for the velocity holds 
		\begin{align*}
			\intTO{ \vu \cdot \partial_{t} \pmb{\varphi} } + 	\intTO{  { \pi(\vr)} \Div \pmb{\varphi} } -\intTO{  \nabla \vu: \nabla \pmb{\varphi} } =- \intO{\vu_0 \cdot \pmb{\varphi}(0,\cdot)}
		\end{align*}
		for any $\pmb{\varphi} \in C_c^1([0,T)\times \T^d) $. 
		\item Moreover, for a.e. $\tau \in (0,T)$ we have the following energy inequality 
		\begin{align*}
			&\int_{\T^3}\lr{	\frac{1}{2} \vert \vu \vert^2+ \Pi(\vr )}(\tau) \;\dx + 	\int_0^\tau \int_{\T^3} \vert \nabla \vu  \vert^2 \; \dx \;\dt 
			\leq 	\int_{\T^3}\lr{	\frac{1}{2} \vert \vu_{0} \vert^2+ \Pi(\vr_{0} )} \;\dx .
		\end{align*}
	\end{itemize}
	
\end{df}


 \subsection{Results for existence of weak solution}
To prove the existence of a solution to problem \eqref{bPF:1}--\eqref{bPF:2}  with the above general form of the pressure \eqref{pre1}--\eqref{pre2}, there is a need of modification such that we will be able to control high oscillation for large values of the density. 
We describe the following modification by introducing the parameter $\mu>0$ which we eventually send to zero
\begin{equation}\label{presser-mu}
	\pi_\mu(\vr) = \pi(\vr) + \mu \vr^\Gamma \mbox{ \ \ for sufficiently large } \Gamma.
\end{equation}
For the pressure $\pi_\mu$, the corresponding pressure potential $\Pi_\mu$ is defined by $\displaystyle  \Pi_\mu(\vr)= \vr \int_{1}^{\vr} \frac{\pi_\mu(\xi)}{\xi^2} \text{d}\xi$. \par 
We also suitably modify the initial data $(\vr_{{\mu}}^0, \vu_{{\mu}}^0)$: 
\begin{align}\label{id-approx1}
	\begin{split}
	&	\vr_{\mu}^0 \in C^2(\T^3)  \text{ satisfies }  0<\underline{\vr} \leq \vr_{\mu}^0(x) \leq \overline{\vr} \text{ for all } x\in \T^3 \text{ and } \vr^{0}_{\mu} \rightarrow \vr_0 \text{ in } L^\gamma(\T^3), \\[4pt]
	& \vu_{\mu}^0 \in C^2(\T^3) \text{ such that } \vu_{\mu}^0 \rightarrow \vu^0 \text{ in }L^2(\T^3) ,\\
	& 	E_{0,\mu}= \int_{\T^3}\lr{	\frac{1}{2} \vert \vu_{\mu}^{0} \vert^2+ \Pi_\mu(\vr^{0}_{\mu} )} \dx  \rightarrow E_0 \text{ as } \mu \rightarrow 0.
	\end{split}
\end{align}

\begin{thm}\label{th:mu}
    Let $\gamma >1$ and $\mu >0$ and let $\vr_\mu^0$ and $\vu_\mu^0$ satisfy \eqref{id-approx1}. 
    For $\Gamma > 3$ the system
    \eqref{bPF:1}--\eqref{bPF:2} with pressure satisfying \eqref{pre1}--\eqref{pre2} and \eqref{presser-mu} admits at least one weak solution following Definition \ref{d:wsol} in the class 
    \begin{align*}
    		&0\leq \vr_\mu \in  C_{\text{weak}}([0,T]; L^{\Gamma}(\T^3)) \text{ \ and \ } \vu_\mu \in  L^2(0,T; W^{1,2}(\T^3) )  \cap C([0,T]; L^{2}(\T^3) ) .
    \end{align*}
\end{thm}
\begin{rmk}
	We note that in Theorem \ref{th:mu}, in order to justify the sense of Definition \ref{d:wsol}, we need to replace $\pi, \Pi$ by $\pi_\mu, \Pi_\mu$ and the initial data $(\vr^0,\vu^0) $ by $(\vr_\mu^0,\vu_\mu^0)$.
\end{rmk}
Theorem \ref{th:mu} is the kernel of this paper, it gives a new construction of weak solutions without parabolic regularization of the continuity equation. It provides the solution to equations with improved pressure. Thus, the next step is to remove the correction $\mu \vr^\Gamma$.

Since the general form of the pressure presented above leads to certain restrictions on the exponent $\gamma$, we also consider another situation with the pressure form
 \begin{align}\label{pr-cmpt}
     &\pi(\vr)= p(\vr)+ q(\vr),\; \text{ with } q\in C_c[0,\infty), \; q(0)=0 \text{ and } \\
     & p(0)=0,\; p(\vr) > 0,\; p^\prime(\vr) >0 \text{ for } \vr>0 \text{ and } p \text{ satisfies  \eqref{pre1}--\eqref{pre2} with }\gamma >1.\nonumber
 \end{align} 

\begin{thm} \label{th:main}
	Let $\gamma > 1 $ and the initial conditions satisfy \eqref{id}. Let the pressure satisfy \eqref{pr-cmpt} with $\gamma>1$ or for the  general pressure satisfying \eqref{pre1}--\eqref{pre2} we have $\gamma\geq \frac{6}{5}$.
	Then there exists $(\vr,\vu)$ a weak solution to \eqref{bPF:1}--\eqref{bPF:2} in the sense of Definition \ref{d:wsol}. 
\end{thm}

In Theorem \ref{th:main}, the hypotheses on $\gamma$ are related to the $L^2$ integrability of the density. They differ from those for the compressible Navier–Stokes equations since, in our case, there is no nonlinear term $\Div(\vr \vu \otimes \vu)$. They correspond to the cases 
$\gamma > 3/2$ and $\gamma \geq 9/5$.

The study of global-in-time solutions to the compressible Navier–Stokes system began relatively late. We should mention the significant impact of Matsumura and Nishida \cite{MN1}, though it was under the restrictive assumption of smallness of the data. An improvement concerning the regularity of the initial data was made in \cite{Dan}, \cite{Mu-Zaj}, \cite{Val-Zaj}.

Since the existence of solutions for small data is not entirely satisfactory, the goal is to define a global-in-time weak solution for the compressible Navier–Stokes equation with large initial data, inspired by the weak solutions for the incompressible case \textit{\`a la Leray}. In this context, it is worth mentioning Hoff's result \cite{H}, which addresses the existence theory for a discontinuous initial data. The existence of weak solutions was later established by Lions \cite{Lions2} and, almost simultaneously, by Feireisl \cite{F}, who expanded the range of $\gamma$. Subsequently, the existence theory for weak solutions was extended to heat-conducting fluids \cite{F}. These techniques have been modified, slightly improved, and applied to various generalizations \cite{FNP}, \cite{FN}, but the core method remains unchanged.

Recently, Bresch and Jabin introduced a new perspective by proposing a method to demonstrate the compactness of the density using the Kolmogorov criterion directly. Essentially, they showed that for an approximate sequence of densities $\{\vr^{(n)}\}_{n \in \mathbb{N}}$, one can establish that
\begin{equation} \label{limit_h} \lim_{|h| \to 0} \sup_{n} \|\vr^{(n)}(\cdot + h,t)-\vr^{(n)}(\cdot,t)\|_{L^1} \to 0. \end{equation} 
This method has found interesting applications in models based on the compressible Navier–Stokes equation (see, e.g., \cite{BJ2}, \cite{BJW}, \cite{BMZ}, \cite{DPSV}, \cite{VaZa}), offering in some cases a clearer understanding of the compactness properties of the studied systems.

However, the first step of the construction, namely the existence of the initial approximate system, remains unsatisfactory. An attempt in this direction was made in \cite{CMZ}, but describing such an approach as `simple' is difficult.

The aim of this paper is to propose a new technique for compactness and the construction of approximate solutions, based on a modification of the Feireisl–Lions technique combined with suitable weights from the Bresch–Jabin approach. The approximate system takes the form described in \cite{CMZ}. We focus our analysis on the simplified system \eqref{bPF:1}–\eqref{bPF:2}, but all the steps can be applied to the full approximation of the compressible Navier–Stokes system, given by the following form: 
\begin{equation} \label{NS_appr} \begin{array}{l} \vr_t + \Div ([\vc{u}]_\epsilon \vr) = 0,\\[4pt] \partial_t(( l + [\vr]_l) \vc{u}) + \Div([\vr \vc{u}]_l \otimes \vc{u}) - \Div \mathbf{S}+ \nabla [p(\vr)]_\epsilon = 0, 
\end{array} \end{equation} where $[f]_l = \kappa_l \ast f$ is the standard mollification by convolution with a smooth positive kernel $\kappa_l$ that tends to the Dirac delta and $\epsilon, l >0$.  We omit the analysis of \eqref{NS_appr}, as the main mathematical challenges related to \eqref{NS_appr} appear to be the same as for the model case \eqref{bPF:1}–\eqref{bPF:2}.

\subsection{The idea of the proof of Theorem \ref{th:mu}}
Large $\Gamma$ (in our case, $\Gamma >3$) removes the problems with insufficient integrability of the density, but also 
gives more fine features for large $\vr$. Namely, in that case, all pressures of type \eqref{presser-mu} can be written in the following very practical way
\begin{equation}
    \pi_\mu(\vr)= \mathsf{q}(\vr) + p_\mu(\vr),\nonumber
\end{equation}
 where $\mathsf{q}$ is compactly supported function and $p_\mu(\vr),\; p_\mu'(\vr),\; p_\mu''(\vr) \geq 0\text{ for } \vr \geq 0 $ with suitable growth conditions.  
 For the pressure $\pi_\mu$, corresponding pressure potential $\Pi_\mu$, defined by $\displaystyle  \Pi_\mu(\vr)= \vr \int_{1}^{\vr} \frac{\pi_\mu(\xi)}{\xi^2} \text{d}\xi$, can be written as
\begin{align}\label{P-mu-0}
	{\Pi}_{\mu}(\vr) =P_\mu(\vr)+\mathsf{Q}(\vr) \text{ with }P_\mu^{\prime \prime \prime}(\vr),\; P_\mu^{\prime \prime }(\vr),\; P_\mu^{\prime }(\vr) \geq 0 \text{ for } \vr \geq 0 \text{ and } \mathsf{Q} \in C_c^1[0,\infty) \cap C^3[0,\infty)  .
\end{align}
The features of the above splitting yields the following important and useful property: there exists $a>0$ such that 
\begin{equation}\label{pP1}
    P_\mu(\vr) \pm a p_\mu(\vr) \mbox{ is convex. \ }
\end{equation}
We will state an exact formulation of $P_{\mu}$ in Lemma \ref{lem:pressure}. We now briefly describe the main steps of the proof of Theorem \ref{th:mu}.
 
\noindent
\begin{itemize}[leftmargin=*]
	\item[-] \textbf{Approximate problem:} We consider the following approximate system
\begin{equation*}
    \begin{array}{l}
         \partial_t \vr + \Div (\vr [\vu]_\ep) + \delta \vr^m=0, \\
         \partial_t \vu-\Delta \vu =-\nabla [\pi_\mu(\vr)]_\ep  ,
    \end{array}
\end{equation*}
where $m>0$ is sufficiently large, $\delta>0$,   $\displaystyle [f]_{{\ep}} = \kappa_\ep \ast f$ and  $\displaystyle [\bf{f}]_{{\ep}} = \kappa_\ep \ast \bf{f}$, where $\kappa_{\ep}$ is the standard mollifier in space as we consider the domain as the case of the torus. The system is accompanied with the initial conditions $(\vr^0_\mu,\vu^0_\mu)$ from \eqref{id-approx1}. In Section \ref{sec3}, we  describe this approximation.

The choice of large $m$ gives us sufficient integrability so that we can proceed with the limit passage $\ep \to 0$. In particular we need to have square integrability of the pressure at this first stage of the approximation. Even in this simplified model the uphill task is to prove the strong convergence of the density, which is crucial to identify the weak limit of the nonlinear composition of the density, in our case the pressure term, i.e., $\displaystyle \lim\limits_{\ep \rightarrow 0} \pi_\mu(\vr_\ep) =\pi_\mu(\vr)$, at least in some weak sense.

	\item[-] \textbf{Continuity and renormalized continuity equation:} With the sufficient regularity of $(\vr,\vu) $ we obtain the following form of the limit continuity equation:
\begin{equation}
    \partial_t \vr + \Div (\vr \vu) + \delta \overline{\vr^m}=0.\nonumber
\end{equation}
Recall that we use the standard notation $\Ov{f(\vr)}$ as a weak limit of sequence $f(\vr_\ep)$ for a weakly convergent sequence $\vr_\ep$ in suitable spaces. 
By the classical method of DiPerna and Lions \cite{DL}, we know that the continuity equation can be viewed in the renormalized sense, i.e., 
\begin{equation}
    \partial_t P_\mu(\vr) + \Div (P_\mu(\vr) \vu) 
    +p_\mu(\vr) \Div \vu+ \delta \overline{\vr^m}
    P'_\mu(\vr)=0. \nonumber
\end{equation}
On the other hand, the direct limit passage in the approximate renormalized continuity equation yields  
\begin{equation}
    \partial_t \overline{P_\mu(\vr)} + \Div (\overline{P_\mu(\vr)} \vu) 
    +\Ov{\Ov{p_\mu(\vr) \Div \vu}}+ \delta \overline{\vr^m
    P'_\mu(\vr)}=0,\nonumber
\end{equation}
where the notation $\Ov{\Ov{f}}$ is introduced to distinguish the possibility of a different limit because of mollification of velocity in the continuity equation.
Therefore, we obtain 
\begin{equation}\label{cont_lim}
    \partial_t (\overline{P_\mu(\vr)}-P_\mu(\vr)) + \Div ((\overline{P_\mu(\vr)} -P_\mu(\vr)) \vu) 
    +\Ov{\Ov{p_\mu(\vr) \Div \vu}}- p_\mu(\vr) \Div \vu+ \delta (\overline{\vr^m
    P'_\mu(\vr)}- \Ov{\vr^m} P'_\mu(\vr)) =0.
\end{equation}
 From the properties of $P_\mu$ as described in \eqref{P-mu-0}, we notice 
 \begin{align*}
 	\overline{P_\mu(\vr)}-P_\mu(\vr) \geq 0 \text{  \ \ and \ \   }\overline{\vr^m
 			P'_\mu(\vr)}- \Ov{\vr^m} P'_\mu(\vr) \geq 0.
 \end{align*}
Therefore, the key term here to analyse is the term
\begin{align}\label{if:kt}
	 \Ov{\Ov{p_\mu(\vr) \Div \vu}}- p_\mu(\vr) \Div \vu .
\end{align}

In this regard we need to understand the role played by an important quantity, the \textit{effective viscous flux}.
\item [-] \textbf{Effective viscous flux:} From the continuity and velocity equation \eqref{bPF:1}--\eqref{bPF:2}, we have 
\begin{align}\label{if:evf1}
	\begin{split}
&	\pt \Div \vu + \Delta \lr{ \pi(\vr) - \Div \vu}=0 \\ 
& \hskip1cm
\text{ which can be stated in the following form } \\
	 &	\pt \lr{\Div \vu - \pi(\vr)} - \Delta{\lr{\Div \vu - \pi(\vr)}}= \Div(\pi(\vr)\vu) + \lr{\vr \pi^\prime(\vr)-\pi(\vr) }\Div \vu. 
	\end{split}
\end{align}

The effective viscous flux $\mathsf{G}$ is given by 
\begin{align}\label{efv}
	\mathsf{G}= \Div \vu - \pi(\vr).
\end{align}
Keeping in mind our approximate system, we need to adapt $\mathsf{G}$ suitably.  For $\ep>0$, we consider 
\[ G_\ep= \Div \vu_\ep -[\pi_\mu(\vr_\ep)]_{\ep}. \]
At least on the formal level, from the equations for  the effective viscous flux \eqref{if:evf1}, we observe that with sufficient regularity of density and velocity, we will be able to obtain compactness (strong convergence) of this quantity $G_\ep$. It plays a very important role in the analysis and it is explained in details in the proof of Lemma \ref{lem:evf:ep}.

\item[-] \textbf{Use of compactness of effective viscous flux in \eqref{if:kt}: } We assume the strong  convergence of $G_\ep$ to $G$. Note that at this level $G\lr{=\Div \vu - \Ov{\pi_\mu(\vr)}}$ might not be equal to $\Div \vu - \pi_\mu(\vr)$. Let us give a closer look at the quantity \eqref{if:kt} and rewrite it as
\begin{align*}
	 & \Ov{\Ov{p_\mu(\vr) \Div \vu}}- p_\mu(\vr) \Div \vu\\
	&=\Ov{\Ov{p_\mu(\vr) p_\mu(\vr)}}- p_\mu(\vr) \Ov{p_\mu(\vr)} + \Ov{\Ov{p_\mu(\vr) \mathsf{q}(\vr)}} - p_\mu(\vr) \Ov{\mathsf{q}(\vr)} +
	(\Ov{p_\mu(\vr)}-p_\mu(\vr))G\\
	&=\!\lr{\Ov{\Ov{p_\mu(\vr) p_\mu(\vr)}}- \Ov{p_\mu(\vr)}^2 } \!\!+ \!{ \lr{\Ov{p_\mu(\vr)}\! - \!p_\mu(\vr)} \Ov{p_\mu(\vr)} }\!+\!\lr{ \Ov{\Ov{p_\mu(\vr) \mathsf{q}(\vr)}} - p_\mu(\vr) \Ov{\mathsf{q}(\vr)} }\!+\!
	\lr{(\Ov{p_\mu(\vr)}-p_\mu(\vr))G}\!.
\end{align*}
The first two terms in the last line of the above equation are non-negative. Also, the last one is of a good form as it depends on the properties of the effective viscous flux $G$. Hence, we need to look closer at the terms with $\mathsf{q}$. So in the weak/distributional sense we have
\begin{align*}
    \lim_{\ep \to 0} \left( [p_\mu(\vr_\ep)]_\ep [\mathsf{q}(\vr_\ep)]_\ep - p_\mu(\vr) \mathsf{q}(\vr_\ep)\right) 
    &=\lim_{\ep \to 0} \left( p_\mu(\vr_\ep)[[\mathsf{q}(\vr_\ep)]_\ep ]_\ep - p_\mu(\vr) [[\mathsf{q}(\vr_\ep)]_\ep]_\ep \right)\\
   & =\lim_{\ep \to 0} \left( (p_\mu(\vr_\ep) - p_\mu(\vr) )[[\mathsf{q}(\vr_\ep)]_\ep ]_\ep \right).
\end{align*}
Since $\mathsf{q}$ is compactly supported, the sequence $\displaystyle [[\mathsf{q}(\vr_\ep)]_\ep ]_\ep $ is uniformly bounded. Also, from \eqref{pP1}, we deduce (see (\ref{id:q-est})) the relation between $p_\mu $ and $P_\mu$ and it yields
\begin{equation*}
    |\Ov{\Ov{p_\mu(\vr) \mathsf{q}(\vr)}} - p_\mu(\vr) \Ov{\mathsf{q}(\vr)}| \leq C (\Ov{P_\mu(\vr)} - P_\mu(\vr)).
\end{equation*}
Thus plugging the expression of $ \displaystyle \lr{\Ov{\Ov{p_\mu(\vr) \Div \vu}}- p_\mu(\vr) \Div \vu}$ into \eqref{cont_lim}, it yields
\begin{equation*}
    \partial_t  (\overline{P_\mu(\vr)}-P_\mu(\vr))
    + \Div ((\overline{P_\mu(\vr)} -P_\mu(\vr)) \vu)  \leq C(\lambda + |G|) (\overline{P_\mu(\vr)}-P_\mu(\vr)) 
\end{equation*}
or  the integrated form 
\begin{align*}
	\Dt \int_{\T^3} \lr{\overline{P_\mu(\vr)}-P_\mu(\vr)} \dx \leq \int_{\T^3} C(\lambda + |G|) \lr{\overline{P_\mu(\vr)}-P_\mu(\vr)} \dx.
\end{align*}
Unfortunately, to use the desired Gr\"onwall's inequality, we need the function $G$ to be point-wise bounded. In our setup, we will be unable to derive the uniform $L^\infty$ bound of the limit effective viscous flux $G$. 
The above inequality seems to be useless. So we need to cheat here a little. 
\item[-] \textbf{Weighted inequality:} We take a weight $w$ such that it is a solution to the problem 
		\begin{align}\label{weight-eqn}
			\pt w + \vu \cdot \nabla w +\Lambda w=0,\; w(0)=1 \mbox{ \ \ with \ \ } \Lambda = 2C(1 + \mathcal{M}(|\nabla \vu|) +|G|).
		\end{align}
First, we note that with the choice of $\Lambda$, the equation for the weight is well-posed. $\mathcal{M}(\cdot)$ stands for the maximal function. Moreover, on the formal level, we derive the following weighted inequality 
\begin{equation}
    \partial_t  \lr{(\overline{P_\mu(\vr)}-P_\mu(\vr)) w }
    + \Div ((\overline{P_\mu(\vr)} -P_\mu(\vr)) \vu w)  \leq \left(C(1 + |G|) - \Lambda \right)(\overline{P_\mu(\vr)}-P_\mu(\vr)) w .\nonumber
\end{equation}
The choice of $\Lambda$ gives us  
\begin{align*}
	\partial_t  \lr{(\overline{P_\mu(\vr)}-P_\mu(\vr)) w }
	+ \Div ((\overline{P_\mu(\vr)} -P_\mu(\vr)) \vu w)  \leq 0
\end{align*}
and eventually 
\begin{equation}
    \int_{\T^3} (\overline{P_\mu(\vr)}-P_\mu(\vr))(\tau) w(\tau) \,\dx =0 \text{ for a.e. } \tau \in (0,T).   \nonumber
\end{equation}
Therefore, we get the identity $ \displaystyle \overline{P_\mu(\vr)}=P_\mu(\vr)$ $w$-a.e. This, however, is not yet sufficient to prove the strong convergence of density.  

\item[-]\textbf{Removal of the weight $w$:} Using the property of $w$ stated in Proposition \ref{prop:weight} we will be able to remove the weight $w$ and as  simple consequence we obtain the strong convergence of $\vr$. But this step requires $L^2$ integrability of the density. We have it at this stage since $\Gamma$ is large. 
 	\item[-] \textbf{Limit passage $\delta \to 0$:} We need to eventually remove the parameter $\delta$ introduced in the continuity equation. In this case, the additional pressure estimate plays a crucial role.  
 \end{itemize}

 \subsection{The idea of the proof of Theorem \ref{th:main}}
 The proof of this theorem is divided into two cases. In the \textit{first case}, we consider a non-monotone pressure of the form \eqref{pr-cmpt}. Meanwhile, the \textit{second case} deals with the pressure \eqref{pre1}--\eqref{pre2} with $\gamma \geq \frac{6}{5}$. 

 \subsubsection*{Case 1: Pressure satisfying \eqref{pr-cmpt}} 
 Here, we follow the idea described in Feireisl \cite{EF2001} for monotone pressure and its subsequent extension for pressure satisfying \eqref{pr-cmpt}, as presented in Feireisl \cite{F2002}.
 \begin{itemize}[leftmargin=*]
     \item[-] To achieve the strong convergence of the density, we consider the \emph{Oscillation defect measure} which is defined as
	\begin{align*}
		\textbf{osc}_\alpha [\vr_\mu \rightarrow \vr ]  := \sup_{k\geq 1} \left(\limsup_{\ep\to0} \int_0^T \int_{\T^3} \vert T_k(\vr_\mu)- T_k(\vr)  \vert^\alpha \,\dx \,\dt  \right),
	\end{align*}
	 $\alpha>0$, and $T_k$ is a suitable truncation function. More precisely, $T_k$ is given by		$T_k(x) = k \;T\left( \frac{x}{k}\right),$
  where $ T $ is a smooth function satisfying
	\begin{align*}
		T(x) = \begin{cases}
			&x,\; 0\leq x\leq 1\\
			& \text{concave smooth},  1\leq x\leq 3\\
			& 2,\; x \geq 3.\\
		\end{cases}
	\end{align*}
	
 \item[-] Next, we establish that for suitable functions $b$, the following effective viscous flux identity holds in the weak sense:
		\begin{align*}
			\overline{\pi(\vr) T_k(\vr) }- \overline{\pi(\vr)}\;\overline{T_k(\vr)} =\overline{T_k(\vr)\dv \vu} - \overline{T_k(\vr)} \dv \vu .
		\end{align*}
 This, combined with the appropriate convergence of approximating sequences, ensures that the oscillation defect is finite for $\alpha=\gamma+1$, i.e.,
  \begin{align}\label{osd-finite}
			\textbf{osc}_{\gamma+1} [\vr_\mu \rightarrow \vr ]<\infty \text{ for } \gamma>1. 
		\end{align}
  \item[-] With the help of \eqref{osd-finite}, we prove the strong convergence of the density, i.e., 
  \[ \Ov{\vr \log \vr}=\vr \log \vr \text{ a.e. in } (0,T)\times \T^3. \]
  Moreover, the renormalized continuity equation also holds in the weak sense.
 \end{itemize}


\subsubsection*{Case 2: Pressure satisfying \eqref{pre1}--\eqref{pre2} with $\gamma \geq \frac{6}{5}$}
In this case, we follow the article of Bresch and Jabin \cite{BJ}.   
\begin{itemize}[leftmargin=*]
    \item[-] The key idea here is to use the Kolmogorov compactness criterion for an approximating sequence of densities. We recall the criterion (Belgacem and Jabin \cite[Lemma 3.1]{BelJ} and Bresch and Jabin \cite[Proposition 1]{BJ2}): 
    \begin{lemma}\label{lem:compact}
	Let $\{X_n\}_{n=1}^\infty$ be a sequence of functions uniformly bounded in $L^p((0,T)\times\T^d)$
	with $1\leq p<+\infty$. Let $ \mathcal{K}_h$ be given by 
 \begin{equation}
  \mathcal{K}_h(z)=
  \left\{ 
  \begin{array}{lr}
   \frac{1}{(|z|+h)^d}& \text{for}\   |z|\leq 1/2, \\[10pt]
   \frac{1}{(1/2 + h)^d} & \text{otherwise}.
  \end{array}
\right.
 \end{equation}
	If $\{\pt X_n\}_{n=1}^\infty$ is uniformly bounded in $L^r(0,T;W^{-1,r}(\T^d))$ with 
	$r\geq 1$ and  
	\eq{\label{criterion}
		\underset{n}{\lim\sup} \left( \frac{1}{\| \mathcal{K}_h\|_{L^1}}  \int_0^T\!\!\!\iintO{
			\mathcal{K}_h(x-y)|X_n(t,x)-X_n(t,y)|^p} \, \dt\right) \to 0, \quad \mbox{ as } h\to 0,
	}
	then, $\{X_n\}_{n=1}^\infty$ is compact in $L^p((0,T)\times\T^d)$.
	Conversely, if $\{X_n\}_{n=1}^\infty$ is compact in $L^p([0,T]\times\T^d)$, 
	then \eqref{criterion} holds.
\end{lemma}
\item[-] Again, we notice that we are unable to establish criterion \eqref{criterion} for an approximate sequence of density $ \{\vr_\mu \}_{\mu>0} $. Instead, with the help of suitable weights $(w_\mu)$ as described in \eqref{weight-eqn} (replacing $\vu$ by $\vu_\mu$), we prove that the approximate solutions to the system satisfy
	\eq{\label{KolCri-rem-w}
		\underset{\mu\to0}{\lim\sup} \left(  \frac{1}{\| \mathcal{K}_h \|_{L^1}}\int_0^T\!\!\!\iintO{
			\mathcal{K}_h(x-y)|\vr_\mu(t,x)-\vr_\mu(t,y)|\; w_\mu(t,x)\;} \, \dt\right) \to 0, \ \mbox{ as } h\to 0.
	}

\item[-] Then, the eventual removal of weight is achieved using suitable properties of $w$ and the square integrability of $\vr$.
 Unfortunately, the restriction $\gamma\geq\frac{6}{5}$ is necessary to ensure it.
\end{itemize}


\section{Proof of Theorem \ref{th:mu}}\label{sec3}

We recall the following approximation steps:
\begin{itemize}[leftmargin=*]
	\item (fixed $\mu $): The aim of Theorem \ref{th:mu} is to prove the existence to the following system
	\begin{align}
		&\pt \vr + \Div(\vr \vu) =0,\label{bPF:1:mu}\\
		&\pt \vu +\nabla {\pi}_{\mu}(\vr) = \Delta \vu\label{bPF:2:mu},
	\end{align}
	where 
	\begin{align}\nonumber
		{\pi}_{\mu}(\vr) ={\pi}(\vr)+ \mu \vr^\Gamma  \text{ with } \Gamma >3. 
	\end{align}
\noindent	

	\item \textbf{($\delta \rightarrow 0$)} The solutions to (\ref{bPF:1:mu})--(\ref{bPF:2:mu}) will be obtained from the following approximate system
	\begin{align}
		&\pt \vr + \Div(\vr \vu)+ \delta \vr^m =0,\label{bPF:1:del}\\
		&\pt \vu +\nabla {\pi}_{\mu}(\vr) = \Delta \vu\label{bPF:2:del}.
	\end{align}

The term $\delta \vr^m$ delivers the basic integrability to control the pressure term in the $L^{2+\ep}$ spaces for some positive $\ep$. Thanks to the additional pressure estimate we will be able to pass with $\delta \to 0$. All this, however, leads to certain restrictions on $m$. On one hand, we need $m$ large to have sufficient regularity of all terms in the momentum equation, on the other hand, in order to pass with $\delta \to 0$ we need to control the density independently of $\delta$ which leads to restrictions from above on $m$. Our technique works well for $m \in (\frac 52 \Gamma +1, \frac 52 \Gamma +2)$.
However, for the sake of clarity, we set $m= \frac 52 \Gamma + \frac 32$.
 
	\item (${\ep} \rightarrow 0$) Regularization of pressure and velocity
	\begin{align}
		&\pt \vr + \Div(\vr [\vu]_{{\ep}})+ \delta \vr^m =0,\label{bPF:1:ep}\\
		&\pt \vu+ \nabla [{\pi}_\mu(\vr) ]_{\ep}= \Delta \vu \label{bPF:2:ep},
	\end{align}
	where $\displaystyle [f]_{{\ep}} = \kappa_\ep \ast f$ and  $\displaystyle [\bf{f}]_{{\ep}} = \kappa_\ep \ast \bf{f}$ with $\kappa_{\ep}$  the standard mollifier in space, we assume  $\kappa_1$ is nonnegative, compactly supported, smooth and even, with the integral equal $1$. In all cases, $\vr(0,\cdot) = \vr^0_\mu$ and $\vu(0,\cdot) = \vu^0_\mu$, see \eqref{id-approx1}. The above mollification makes our problem easier from the viewpoint of the construction of solutions. Then the main task stays the limit $\ep \to 0$.

\end{itemize}

\subsection{Existence for fixed $\ep$, $\delta$, $\mu$}

We consider the initial data $(\vr_{{\ep},\delta,{\mu}}^0, \vu_{{\ep},\delta,{\mu}}^0)$ with $\vr_{\ep,\delta,{\mu}}^0= \vr_{\mu}^0$, $\vu_{\ep,\delta,{\mu}}^0= \vu_{\mu}^0$  such that they satisfy \eqref{id-approx1}. 
\begin{thm}\label{exi-fix}
	Let $\ep,\,\delta,\, $ and $\mu$ be fixed. Then there exists a strong solution $(\vr,\vu)$ to system \eqref{bPF:1:ep}--\eqref{bPF:2:ep} with initial conditions fulfilling \eqref{id-approx1} such that
	\begin{equation*}
		\vu \in W^{1,p}(0,T; L^{p}(\T^3))\cap L^p(0,T;W^{2,p}(\T^3)), \qquad {0\leq}\,\vr \in W^{1,p}(0,T; L^{p}(\T^3))\cap L^p(0,T;W^{1,p}(\T^3)) 
	\end{equation*}
	for any $p<\infty$. 
\end{thm}

For all parameters fixed we present the proof of the theorem above in Appendix in Section \ref{App.1}.

 \subsubsection{Modified pressure and its properties }

Let us recall the pressure $
\pi_\mu$ and the pressure potential $(\Pi_\mu)$: 
\begin{align}
	{\pi}_{\mu}(\vr) =\pi(\vr)+ \mu \vr^\Gamma \text{ \ \ and \ \ } {\Pi}_{\mu}(\vr) =\Pi(\vr)+ \frac{\mu}{\Gamma-1} \vr^\Gamma  . \nonumber
\end{align}

\begin{lemma}\label{lem:pressure}
Given $\pi_\mu$ and $\Pi_\mu$, there exist $P_\mu$ and $\mathsf{Q}$ such that ${\Pi}_{\mu}(\vr) =P_\mu(\vr)+\mathsf{Q}(\vr) \text{ with }$
\begin{align}\label{P-mu}
	P_\mu^{\prime \prime \prime}(\vr),\; P_\mu^{\prime \prime }(\vr),\; P_\mu^{\prime }(\vr),\; P_\mu(\vr) \geq 0 \text{ for } \vr \geq 0 \text{ and } \mathsf{Q} \in C_c^1[0,\infty) \cap C^3_c[0,\infty)  .
\end{align}
Moreover, we define the corresponding pressures $\pi_\mu(\vr)= p_\mu(\vr) + \mathsf{q(\vr)}$ as 
\begin{align}\label{rf-p-mu}
	p_\mu(\vr) = \vr P_\mu^\prime(\vr)-P_\mu(\vr) \text{ and } \mathsf{q}(\vr)= \vr \mathsf{Q}^\prime(\vr)-\mathsf{Q}(\vr). 
\end{align} 
As a consequence, we have the following properties of $P_\mu$ and $\mathsf{Q}$:
\begin{itemize}[leftmargin=*]
	\item There exists $C_1,\; C_2, \; C_3 > 0$, such that 
	\begin{align*}
		&C_1 \vr^\Gamma	\leq P_\mu(\vr)  \leq C_2 \lr{ \vr^\Gamma+ \vr^\gamma+1} \text{ for } \vr \geq 0, \\
		&	C_1 \vr^{\Gamma-1} \leq P_\mu^\prime(\vr)  \leq C_2 \lr{ \vr^{\Gamma-1}+ \vr^{\gamma-1}+1} \text{ for } \vr \geq 0.
	\end{align*}
	Furthermore, we obtain 
	\begin{align*}
		0\leq p_\mu (\vr)  \leq C_3 \lr{ \vr^\Gamma+ \vr^\gamma+1} \text{ and } 0\leq p_\mu^\prime(\vr)    \leq C_3 \lr{ \vr^{\Gamma-1}+ \vr^{\gamma-1}+1} \text{ for } \vr \geq 0.
	\end{align*}
	\item For $\mathsf{Q} $ and $\mathsf{q}$ it holds 
		$|\mathsf{Q}(\vr) |  + |\mathsf{Q}^{\prime }(\vr) |  \leq  C \text{ for } \vr \geq 0$.
	\item Moreover, there exists $\lambda_{\mathsf{q}} >0 $ such that 
	\begin{align}\label{P-p-mu}
		\lambda_{\mathsf{q}}  P_\mu(\vr) \pm p_\mu(\vr) \text{ is convex.}
	\end{align}
	
\end{itemize}
\end{lemma}
 
We skip the proof of the above lemma. It is based on an easy observation that for some large $R> 1$, it holds
	$\partial_\vr^{(i)}\lr{P(\vr) +\frac{\mu}{\Gamma-1} \vr^\Gamma } >0$  for $ \vr >R$, $ i=0,1,2,3.$
	Hence, we are able to write the pressure potential $\Pi_\mu$ such that it satisfies \eqref{P-mu}. 
	There exists an $M>R>0$ such that we write 
	\begin{align*}
		P_\mu(\vr)= \chi(\vr) P(\vr) + \frac{\mu}{\Gamma-1} \vr^\Gamma ,
	\end{align*}
	where 		$\chi \in C^3[0,\infty), \; \chi^\prime \geq 0$ 
 and $\chi$ is zero for $\vr < M$
 and equal $1$ for $\vr > 2M$, potential $P_\mu$ fulfills the lemma.


\subsection{Limit passage in step $ \ep \rightarrow 0 $ }
For fixed $\delta, \mu$, we denote $(\vr_{{\ep},\delta,\mu},\vu_{{\ep},\delta,\mu})=(\vr_{\ep},\vu_{\ep})$. We will now perform the limit passage $\ep \rightarrow 0 $ and prove the following theorem: 

\begin{lemma} 
	Let $\delta$ and ${\mu}$ be fixed. Then for initial data satisfying \eqref{id-approx1} there exists unique $(\vr_{\delta,\mu},\vu_{\delta,\mu})$ to system \eqref{bPF:1:del}--\eqref{bPF:2:del} in the class
	\begin{align*}
		&0\leq \vr_{\delta,{\mu}} \in  C_{\text{weak}}([0,T]; L^{p_1\Gamma}(\T^3)) \cap L^{p_2\Gamma} ((0,T)\times \T^3) \text{ with } p_1 =\frac{m-1}{\Gamma} > \frac{5}{2} \text{ and } p_2=\frac{2(m-1)}{\Gamma} > 5\\ 
		&   \text{ and } \vu_{\delta,{\mu}} \in  L^2(0,T; W^{1,2}(\T^3) )  \cap C([0,T]; L^{2}(\T^3) ) 
	\end{align*}
	such that 
	\begin{align*}
		&\intTO{ \Big(\vr_{\delta,{\mu}} \partial_t \psi + \vr_{\delta,{\mu}} \vu_{\delta,{\mu}} \cdot \nabla \psi\Big)} -\intTO{ \delta { \vr_{\delta,{\mu}}^m} \psi}= -\intO{\vr^{0}_\mu \psi(0)} 
	\end{align*}
	for any $\psi \in  C_c^1([0,T)\times \T^3) $ and 
	\begin{align*}
		\intTO{ \vu_{\delta,{\mu}} \cdot \partial_{t} \pmb{\varphi} } + 	\intTO{  { {\pi}_\mu(\vr_{\delta,{\mu}})} \Div \pmb{\varphi} } -\intTO{  \nabla \vu_{\delta,{\mu}} : \nabla \pmb{\varphi} } =- \intO{\vu^{0}_{\mu} \cdot\pmb{\varphi}(0,\cdot)}
	\end{align*}
	for any $\pmb{\varphi} \in C_c^1([0,T)\times \T^3) $. 
	Moreover, for a.e. $\tau \in (0,T)$ we have the following energy inequality 
	\begin{align*}
		&\int_{\T^3}\lr{	\frac{1}{2} \vert \vu_{\delta,{\mu}} \vert^2+ \Pi_{{\mu}}(\vr_{\delta,{\mu}} )}(\tau) \;\dx + 	\int_0^\tau \int_{\T^3} \vert \nabla \vu_{\delta,{\mu}}  \vert^2 \; \dx \;\dt+\int_0^\tau \int_{\T^3}  \delta  \vr_n^{m} \Pi_{{\mu}}^\prime(\vr_{\delta,\mu}) \;\dx\; \dt \nonumber \\
		&\leq 	\int_{\T^3}\lr{	\frac{1}{2} \vert \vu^{0}_{\mu} \vert^2+ \Pi_{\mu}(\vr^{0}_{\mu} )} \;\dx 
	\end{align*}
	and the renormalized continuity equation
	\begin{align*}
		&\intTO{ \Big(b(\vr_{\delta,{\mu}}) \partial_t \psi + b(\vr_{\delta,{\mu}}) \vu_{\delta,{\mu}} \cdot \nabla \psi\Big)} + \intTO{(b'(\vr_{\delta,{\mu}}) \vr_{\delta,{\mu}}-b(\vr_{\delta,{\mu}})) \Div \vu_{\delta,{\mu}} \psi} \\
  &-\intTO{ \delta { \vr_{\delta,{\mu}}^m} b'(\vr_{\delta,{\mu}}) \psi}= -\intO{b(\vr^{0}_\mu) \psi(0)} 
	\end{align*}
	holds for any $b \in C^1([0,\infty))$ such that $b'(z) = 0$ for $z \geq M_0$ for a certain $M_0>0$.
\end{lemma}

\subsubsection{Uniform estimates} 

The first step is to derive the uniform estimates with respect to $\ep$. In this regard,  we first consider $\Gamma>3$ with the choice  
\begin{align}
\label{m-Gamma}
	 m =  \frac{5}{2}\Gamma +\frac 32 .
\end{align}

\begin{lemma}\label{lem:ep:uniest}
	Let the initial data follow \eqref{id-approx1} and let $m$ and $\Gamma$ be related by \eqref{m-Gamma}, then the solution $(\vr_\ep,\vu_\ep)$ and the effective viscous flux $G_{\ep}\lr{= \Div \vu_{\ep} - [p_{{\mu}}(\vr_{\ep})]_{\ep}}$ satisfy the following estimates: 
	\begin{align}\label{unif-est-ep}
		\begin{split}
		&\Vert  \vu_{\ep}  \Vert_{L^2_t W^{1,2}_x}+ \Vert \vu_{\ep}  \Vert_{L^\infty_t L^{2}_x}+ \Vert   \vr_{\ep}^\Gamma  \Vert_{L^\infty_t L^{1}_x} \leq C\\
		& \Vert \pi_{\mu}(\vr_{\ep}) \Vert_{L^\infty_t L^{p_1}_x} + 	\Vert \pi_{\mu}(\vr_{\ep}) \Vert_{ L^{p_2}_{t,x}} + \Vert    \vr_{\ep}^{m } \pi_{\mu}^\prime(\vr_{\ep})\Vert_{L^s_{t,x}}  \leq C \\
		&	\Vert \partial_{t} (p_{\mu}(\vr_{\ep})) \Vert_{L^s_t W^{-1,s}_x} + \Vert \pt \vu_{\ep}  \Vert_{L^2_t W^{-1,2}_x} +	\Vert 	\pt (-\Delta)^{-1} \Div  \vu_{\ep} \Vert_{L^2_{t,x}}  \leq C\\
		& 	\Vert  G_{\ep} \Vert_{L^2_{t,x}} \leq C,
		\end{split}
	\end{align}
	where $p_1, p_2$ and $s$ satisfy 
	\begin{align}\label{p-s}
		p_1=\frac{m-1}{\Gamma} > \frac{5}{2},\; p_2=\frac{2(m-1)}{\Gamma} > 5 \text{ \ and \ } s= \frac{2(m-1)}{m+\Gamma-1} >\frac{10}{7}.
	\end{align}
	\end{lemma}
 
	\textbf{Proof of Lemma \ref{lem:ep:uniest}: }
First, from the approximate continuity equation we have 
\begin{align*}
	\frac{d}{dt}	\int_{\T^3} \vr_{\ep} \;\dx +  \int_{\T^3} \delta  \vr_{\ep}^{m} \;\dx = 0 
	\text{ \ \ and it yields \ \ }
	\Vert \vr_{\ep} \Vert_{L^\infty_t L^1_x} + \Vert \vr_{\ep} \Vert_{L^m_{t,x} } \leq C.
\end{align*}
We write the renormalized continuity equation by multiplying the approximate continuity equation by $b^\prime(\vr_{\ep}) $. It reads
\begin{align}
	\pt b(\vr_{\ep}) + \Div(b(\vr_{\ep}) [\vu_{\ep}]_{\ep})+ \lr{\vr_{\ep} b^\prime(\vr_{\ep})-b(\vr_{\ep})}[ \Div \vu_{\ep}]_{\ep}+  \delta \vr_{\ep}^{m} b^\prime(\vr_{\ep}) =0,\nonumber
\end{align}
where $b$ is sufficiently smooth, whose derivative has bounded support. By approximation, we may extend its validity to less smooth functions with a controlled behaviour at zero and infinity. 
Using the property $\displaystyle \vr \Pi_{\mu}^\prime (\vr) -\Pi_{{\mu}}(\vr)=\pi_{\mu}(\vr )$, we take $b(\vr_\ep)= \Pi_\mu(\vr_\ep)$ 
\begin{align}\label{P-n1}
	\pt \Pi_{\mu}(\vr_{\ep}) + \Div(\Pi_{\mu}(\vr_{\ep}) [ \vu_{\ep}]_{\ep})+  \pi_{{\mu}}(\vr_{\ep})[ \Div \vu_{\ep}]_{\ep}+  \delta  \vr_{\ep}^{m} \Pi_\mu^\prime(\vr_{\ep})=0. 
\end{align}
Next, we observe that 
\begin{align*}
	\int_{\T^3} \pi_{{\mu}}(\vr_{\ep})[ \Div \vu_{\ep}]_{\ep} \dx = - \int_{\T^3} \vu_\ep \cdot \nabla  [\pi_{{\mu}}(\vr_{\ep})]_{\ep}  \,\dx. 
\end{align*}
Therefore, from the momentum equation \eqref{bPF:2:ep} and the equation for pressure potential \eqref{P-n1}, we get
\begin{align*}
	\frac{1}{2} \frac{d}{dt}\int_{\T^3}	\vert \vu_{\ep} \vert^2 \;\dx +   \frac{d}{dt}	\int_{\T^3} \Pi_{\mu}(\vr_{\ep}) \;\dx + 	\int_{\T^3} \vert \nabla \vu_{\ep} \vert^2 \; \dx + \int_{\T^3}  \delta  \vr_{\ep}^{m} \lr{ \Pi_{{\mu}}^\prime(\vr_{\ep})+C} \;\dx \leq C_1 \int_{\T^3} \delta  \vr_{\ep}^{m}\; \dx ,
\end{align*}
where $ \Pi_{{\mu}}^\prime(\vr_{\ep}) \geq C \vr_{\ep}^{\Gamma-1}$. It yields
\begin{align}\label{bound-n1}
	\Vert  \vu_{\ep}  \Vert_{L^2_t W^{1,2}_x}+ \Vert \vu_{\ep}  \Vert_{L^\infty_t L^{2}_x}+ \Vert   \vr_{\ep}^\gamma  \Vert_{L^\infty_t L^{1}_x} +{\mu}  \Vert   \vr_{\ep}^\Gamma  \Vert_{L^\infty_t L^{1}_x} \leq C.
\end{align}
Considering the renormalized continuity equation with $b = \vr^\Gamma$ (again, this choice is justified) we have
\begin{align*}
	\pt \vr_{\ep}^\Gamma + \Div( \vr_{\ep}^\Gamma  [\vu_{\ep}]_{\ep})+ (\Gamma-1)\vr_{\ep}^\Gamma [ \Div \vu_{\ep}]_{\ep}+  \delta \Gamma \vr_{\ep}^{m+\Gamma-1} =0. 
\end{align*}
Moreover, for any $m-1\geq q\geq \Gamma  $ the following equation holds:
\begin{align*}
	\pt \vr_{\ep}^q + \Div(\vr_{\ep}^q [ \vu_{\ep}]_{\ep})+ (q-1) \vr_{\ep}^q [\Div \vu_{\ep}]_{\ep}+  \delta q   \vr_{\ep}^{m+q -1 }=0.
\end{align*}
Eventually, using the fact $ 2 q \leq q+m-1$ and $\displaystyle \int_{\T^3}\vert \vr_{\ep}^q [\Div \vu_{\ep}]_{\ep} \vert \dx \leq \Vert  \kappa_{\ep} \Vert_{L^1} \Vert  \vr_{\ep}^q \Vert_{L^2_x} \Vert  \Div \vu_{\ep} \Vert_{L^2_x}$, we get 
\begin{align*}
	\frac{d}{dt} \intO{\vr_{\ep}^q} + \frac{\delta q}{2} \intO{ \vr_{\ep}^{m+q -1 } } \leq 1+ c(\delta,q) \intO{\vert \nabla \vu_{\ep} \vert^2}.
\end{align*}
Taking $q=m-1 $, with the help of \eqref{m-Gamma}, we obtain ($2m-2> 5\Gamma)$
\begin{align}\label{bound-n2}
	&	\Vert  \vr_{\ep}^\Gamma \Vert_{L^\infty_t L^{p_1}_x} + 	\Vert  \vr_{\ep}^\Gamma \Vert_{ L^{p_2}_{t,x}} + \Vert    \vr_{\ep}^{m+\Gamma -1 } \Vert_{L^s_{t,x}}  \leq C.
\end{align}
We recall $\displaystyle \pi_{\mu}(\vr) = \pi(\vr) + {\mu} \vr^\Gamma $, where $\pi(\vr) $ satisfies \eqref{pre1}--\eqref{pre2}, hence by Lemma \ref{lem:pressure}
\begin{align*}
	\Vert \pi_{\mu}(\vr_{\ep}) \Vert_{L^\infty_t L^{p_1}_x} + 	\Vert \pi_{\mu}(\vr_{\ep}) \Vert_{ L^{p_2}_{t,x}} + \Vert    \vr_{\ep}^{m } \pi_{\mu}^\prime(\vr_n)\Vert_{L^s_{t,x}}  \leq C .
\end{align*} 
We write the equation for $p_\mu$ as defined in \eqref{rf-p-mu}
\begin{align*}
	\pt p_{\mu}(\vr_{\ep}) + \Div(p_{\mu}(\vr_{\ep}) [\vu_{\ep}]_{{\ep}})+ \lr{\vr_{\ep}p_{\mu}^\prime(\vr_{\ep})-p_{\mu}(\vr_{\ep})} [\Div \vu_{\ep}]_{\ep}+  \delta \vr_{\ep}^{m}p_{\mu}^\prime(\vr_{\ep}) =0. 
\end{align*}
Our next goal is to derive an estimate of $\partial_{t} (p_{{\mu}} (\vr_{\ep}))$. To achieve this, we notice the following estimates: 
\begin{align*}
	\Vert  \vr_{\ep}^{m } p_{\mu}^\prime(\vr_\ep) \Vert_{L^s_{t,x}} +
	\Vert  \lr{\vr_{\ep}p_{\mu}^\prime(\vr_{\ep})-p_{\mu}(\vr_{\ep})} [ \Div \vu_{\ep}]_{\ep} \Vert_{L^{s}_{t,x}} +
	\Vert p_{\mu}(\vr_{\ep}) [\vu_{\ep}]_{\ep} \Vert_{L^\frac{2(m-1)}{\Gamma}_{t}L^{s}_{x} \cap L^s_{t}L^{\frac{6(m-1)}{m-1+3\Gamma }}_{x}} \leq C,
\end{align*}
where we  use the structural assumption of $p$.
Hence, from the above estimate along with \eqref{bound-n1} we have 
\begin{align}\label{bound-n3}
	\Vert \partial_{t} (p_{\mu}(\vr_{\ep})) \Vert_{L^s_t W^{-1,s}_x} + \Vert \pt \vu_{\ep}  \Vert_{L^2_t W^{-1,2}_x} \leq C.
\end{align} 
The effective viscous flux is on this level of approximation $G_{\ep}= \Div \vu_{\ep} - [p_{{\mu}}(\vr_{\ep})]_{\ep}$. It satisfies
\begin{align}\label{evf2}
	\pt  (-\Delta)^{-1} \Div \vu_{\ep} = -G_{\ep}.
\end{align}
From \eqref{bound-n1} and \eqref{bound-n2} we have 
\begin{align}\label{bound-n4}
	\Vert  G_{\ep} \Vert_{L^2_{t,x}} \leq C, \ 
	\text{therefore by \eqref{evf2} } 
	\Vert 	\pt (-\Delta)^{-1} \Div  \vu_{\ep} \Vert_{L^2_{t,x}} \leq C.
\end{align}
Collecting all the above estimates we complete the proof of Lemma \ref{lem:ep:uniest}.

\subsubsection{Convergence}\label{subsec:ep:conv}
With the help of Lions--Aubin lemma and \eqref{unif-est-ep} proved in Lemma \ref{lem:ep:uniest}, we have following convergences: 
\begin{align}\label{conv-n1}
	\begin{split}
		& \vu_{\ep} \rightharpoonup \vu  \text{ weak-{*}ly in }  L^\infty_t L^2_x,\\
		& \vu_{\ep} \rightharpoonup \vu  \text{ weakly in }  L^2_t W^{1,2}_x, \\
		& \pt \vu_{\ep} \rightharpoonup \pt \vu  \text{ weakly in }  L^2_t W^{-1,2}_x, \\
		& \vu_{\ep} \rightarrow \vu  \text{ strongly in }  L^2_t L^{2}_x.
	\end{split}
\end{align}
Also, we deduce that 
\begin{align}\label{conv-n2}
	\begin{split}
		&\Pi_{\mu}(\vr_{\ep}) \rightharpoonup\Ov{\Pi_{\mu}(\vr) } \text{ weak-{*}ly in } L^{\infty}_t L^{p_1}_{x} \cap L^{p_2}_{t,x},\\
		&\pi_{{\mu}}(\vr_{\ep}) \rightharpoonup \Ov{\pi_{\mu}(\vr) } \text{ weak-{*}ly in } L^{\infty}_t L^{p_1}_{x} \cap L^{p_2}_{t,x}.
	\end{split}
\end{align}
In particular, we use the following notation:
\begin{align*}
	&p_{{\mu}}(\vr_{\ep}) \rightharpoonup \Ov{p_{\mu}(\vr) } \text{ weak-{*}ly in } L^{\infty}_t L^{p_1}_{x} \cap L^{p_2}_{t,x},\\
	&\mathsf{q}(\vr_{\ep}) \rightharpoonup \Ov{\mathsf{q}(\vr) } \text{ weak-{*}ly in } L^{\infty}_{t,x}.
\end{align*}
Next, we note that for $f_\ep$ with $ f_{\ep} \rightarrow f  \text{ weakly/weak-{*}ly in } L^p_{t,x}$ and $p>1$ as $\ep \rightarrow 0$, we have  
\begin{align*}
	\int_0^T\int_{\T^3} [f_\ep]_\ep \psi\; \dx;\dt  = 	\int_0^T \int_{\T^3}  f_\ep [\psi]_\ep \;\dx; \dt \rightarrow \int_0^T\int_{\T^3}  f \psi \;\dx 
\end{align*}
for any $ \psi \in L^{p^\prime}_{t,x}$ as $\ep \rightarrow 0$. This follows from the fact that $[ \psi]_{\ep} \rightarrow \psi  \text{ strongly in } L^{p^\prime}_{t,x}$. Therefore, we obtain 
\begin{align}
	\begin{split}
		&[p_{{\mu}}(\vr_{\ep})]_{\ep} \rightharpoonup \Ov{p_{\mu}(\vr) } \text{ weakly in } L^{p_2}_{t,x},\\
		&G_\ep \rightharpoonup  G(= \Div \vu -\Ov{\pi_{\mu}(\vr)})  \text{ weakly in }  L^2_t L^{2}_x,\\
		&[\vu_{\ep}]_{\ep} \rightharpoonup \vu  \text{ weakly in }  L^2_t W^{1,2}_x \text{ and }\text{ strongly in }  L^2_t L^{2}_x. \nonumber
	\end{split}
\end{align}
Moreover, we denote
\begin{align*}
	\lr{[p_{\mu}(\vr_\ep)]_\ep }^2\rightharpoonup \Ov{\Ov{p_\mu(\vr)^2}}=:\mathcal{P} \text{ weakly in } L^{p_1}_{t,x}.
\end{align*}
Since $F(\lambda) = \lambda^2 $ is convex, we obtain the following inequality in the sense of distributions:
\begin{align*}
	\mathcal{P} \geq \lr{  \Ov{p_{\mu}(\vr) }   }^2.
\end{align*}


\subsubsection{First limit passage}
\begin{itemize}[leftmargin=*]
	\item \textbf{Continuity equation}: After the limit passage in the weak form of \eqref{bPF:1:ep},  we have
	\begin{align}\label{38a}
		&\intTO{\Big(\vr \partial_t \psi + \vr \vu \cdot \nabla \psi\Big)} -\intTO{ \delta \Ov{ \vr^m} \psi}= -\intO{\vr^{0}_\mu \psi(0)}
	\end{align}
 \text{ for all } $\psi \in C_c^1([0,T)\times \T^3) $.
	\item \textbf{Momentum equation} after limit passage in the weak form of \eqref{bPF:2:ep} reads as
	\begin{align*}
		\intTO{ \vu \cdot \partial_{t} \pmb{\varphi} } + 	\intTO{  \Ov{ \pi_{\mu}(\vr)} \Div \pmb{\varphi} } -\intTO{  \nabla \vu : \nabla \pmb{\varphi} } =- \intO{\vu^0_\mu \cdot\pmb{\varphi}(0,\cdot)}
	\end{align*}
	\text{ for all } $\pmb{\varphi} \in C_c^1([0,T)\times \T^3) $.
	\item Letting $\ep \rightarrow 0$ in the equation for pressure potential \eqref{P-n1}, we obtain (it is rather equation for $P_\mu$ not $\Pi_\mu$, but we can just replace $\pi_\mu$ by $p_\mu$ and we are in more or less in the same situation)
	\begin{align}\label{P2}
		\pt \Ov{P_{{\mu}}(\vr)} + \Div(\Ov{P_{\mu}(\vr)} \vu)+  \Ov{\Ov{p_{{\mu}}(\vr) \Div \vu}}+  \delta\gamma \Ov{\vr^{m}P_{{\mu}}^\prime(\vr)} =0,
	\end{align}
	where 
	\begin{align*}
		p_{{\mu}}(\vr_{\ep})[ \Div \vu_{\ep}]_{\ep} \rightharpoonup\Ov{\Ov{p_{{\mu}}(\vr) \Div \vu}} \text{ weakly in } L^{p}_{t,x} \text{ with some } p>1.
	\end{align*}
	We use different notation because $ \Ov{\Ov{p_{{\mu}}(\vr) \Div \vu}}  $ may not be equal to $\Ov{p_{{\mu}}(\vr) \Div \vu}$. 

	We further have $ \displaystyle \lim\limits_{\ep\rightarrow 0} \int_0^T\int_{\T^3} \lr{	[p_{{\mu}}(\vr_{\ep}) \psi]_{\ep } - [p_{{\mu}}(\vr_{\ep})]_{\ep } \psi} \Div \vu_{\ep}\,\dx \,\dt
	=0 $, which follows from a simple version of Friedrich's commutator lemma.

	\item Due to \eqref{38a} and as $\vr \in L^2_{t,x}$, we get the renormalized form of \eqref{38a} in the form
	\begin{align*}
		&\pt b(\vr) + \Div(b(\vr) \vu)+ (\vr b^\prime(\vr)-b(\vr))\Div \vu + \delta \Ov{\vr^m} b^\prime(\vr)=0, 
	\end{align*}
 where the equation is to be understood in the weak sense. 
	Considering $b= P_{{\mu}}$, we obtain 
	\begin{align}\label{P1}
		&\pt P_{{\mu}}(\vr) + \Div(P_{{\mu}}(\vr) \vu)+p_{\mu}(\vr )\Div \vu + \delta \gamma \Ov{\vr^m} P_{{\mu}}^\prime(\vr)=0. 
	\end{align}
	
\end{itemize}


\subsubsection{Strong convergence of density}
First we establish that 
\begin{align}\label{c-w-1}
	\intO{	 \lr{\Ov{P_{\mu}(\vr)} -P_{\mu}(\vr)} (\tau) w(\tau) }  =  0 \text{ for a.e. } \tau \in (0,T)
\end{align} 
for some suitable weight function $w$ and then we eventually  remove the weight and derive the desired identity 
\begin{align*}
	\intO{	 \lr{\Ov{P_{\mu}(\vr)} -P_{\mu}(\vr)} (\tau)  }  =  0 \text{ for a.e. } \tau \in (0,T)
\end{align*} 
which in turn yields the strong convergence of the density.\par

The particular weight function we are considering can be described by the following proposition. 
\begin{prop}\label{prop:weight}
	Given $\displaystyle \vu \in L^2(0,T; W^{1,2}(\T^3;\R^3)$ and, $\Lambda= a_1+a_2\mathcal{M}(|\nabla \vu|)+ a_3 |G|$ with $ 0\leq a_1,a_2,a_3$,\;  $G\in L^2((0,T) \times \T^3)$ and $\vr\in L^2((0,T) \times \T^3)$, then there exists a function $w$ which satisfies  
	\begin{align}\nonumber
		\pt w + \vu \cdot \nabla w +\Lambda w=0,\; w(0)=1.
	\end{align}
	Moreover, the following properties hold:
	\begin{itemize}
		\item[1.]  For any $(t,x)\in(0,T)\times \T^3$, $0\leq w(t,x) \leq 1$.
		\item[2.] Furthermore, we have
		\begin{equation}\label{boundlogw}
			\underset{t>0}{{\rm ess}\sup} \intO{ \vr |\log w| } \leq C  .
		\end{equation}
	\end{itemize}
	\end{prop}
 \textit{Proof.} Here we proceed with the following steps: 
 \begin{itemize}[leftmargin=*]
     \item At first for general $\vu \in L^2(0,T; W^{1,2}(\T^3;\R^3))$ and $0< \Lambda \in L^2((0,T) \times \T^3)$, we consider a space-time regularization of $\vu $ with standard mollifier and denote it by $ [\vu ]_\ep $. Moreover, we consider $\Lambda_\ep= a_1+ a_2\mathcal{M}(|\nabla [\vu]_\ep | + a_3 [|G|]_\ep) $. Using the method of characteristics, for each  $\ep>0$, we have existence of smooth solution 
     \begin{align*}
        \pt w_\ep + [\vu]_\ep \cdot \nabla w_\ep +\Lambda_\ep w_\ep=0,\;  w_\ep|_{t=0}=1.
     \end{align*}
     \item We observe that $w_\ep \geq 0$ and for any $1\leq p<\infty$, we have 
     \begin{align*}
         \frac{d}{dt} \int_{\T^3} |w_\ep|^p \dx +   \int_{\T^3}\lr{-\frac{1}{p}\Div [\vu]_\ep + \Lambda_\ep} |w_\ep|^p  \dx =0. 
     \end{align*}
     Therefore, $w_\ep\leq 1$. 
     \item The weak limit of $w_\ep$ is denoted by $w$, and using standard commutator estimates, we conclude the existence of the weak solution of \eqref{weight-eqn}.
     \item For the second part, on the formal level the equation for $\vr|\log w|$ is equal to
		\begin{align}\label{eq:rhologw}
			\pt\lr{\vr|\log w|}+\Div(\vr|\log w|\vu)+\delta \vr^m|\log w|= \vr \Lambda.
		\end{align}
		This yields
		\begin{align*}
			\frac{d}{dt} \intO{\vr|\log w|} + \delta \intO{\vr^m|\log w| } \leq \intO{  \vr \Lambda}.
		\end{align*}
 \end{itemize}



\subsubsection*{Step 1: Proof of weighted identity \eqref{c-w-1}} 
From \eqref{P1} and \eqref{P2},
for $\psi \in C^1 ([0,\tau]; C^1(\T^d)$ we have 
\begin{align}\label{wfbarp-p}
	&\intO{	 \lr{\Ov{P_{\mu}(\vr)} -P_{\mu}(\vr)} (\tau) \psi(\tau) } - \inttauO{	\lr{\Ov{P_{\mu}(\vr)} -P_{\mu}(\vr)} \pt \psi }\nonumber\\
	&- \inttauO{ \lr{\Ov{P_{\mu}(\vr)}-P_{\mu}(\vr)}  \vu \cdot \nabla \psi }+\inttauO{\lr{\Ov{\Ov{p_{\mu}(\vr) \Div \vu}}-p_{\mu}(\vr) \Div \vu} \psi } \nonumber \\
	&+  \inttauO{\delta  \lr{\Ov{ \vr^{m}P_{\mu}^\prime(\vr)} {-\Ov{\vr^m}P_{\mu}^\prime(\vr)} }\psi }	= \intO{	 \lr{\Ov{P_\mu(\vr)}(0,x) -P_\mu(\vr^0_\mu)} \psi(0,x)}=0. 
\end{align}
Next, we  replace the term $ \displaystyle \inttauO{\lr{\Ov{\Ov{p_{\mu}(\vr) \Div \vu}}-p_{\mu}(\vr) \Div \vu} \psi}   $ using the next lemma: 
\begin{lemma}\label{lem:wi:ep}
	For any $0<\tau <T$ and $\psi \in C^1([0,\tau]; C^1(\T^3))$ we have the following identity: 
	\begin{align}\label{id:pdiv}
		&\inttauO{\lr{\Ov{\Ov{p_{\mu}(\vr) \Div \vu}}-p_{\mu}(\vr) \Div \vu} \psi } \nonumber \\
		&= 	\inttauO{G \lr{\Ov{p_{\mu}(\vr)} -p_{\mu}(\vr) }  \psi } + 	\inttauO{ \lr{\mathcal P
				- \Ov{p_{\mu}(\vr)} p_{\mu}(\vr)}\psi  } \\
		&\quad +\inttauO{ \lr{\Ov{\Ov{\mathsf{q}(\vr) p_\mu(\vr)}}  -\Ov{\mathsf{q}(\vr)} p_\mu(\vr)} \psi } \nonumber. 
	\end{align}
\end{lemma}
One of the key ingredients we use here is the compactness of the effective viscous flux. More precisely, the following lemma along with the convergences mentioned in Subsection \ref{subsec:ep:conv} is enough for to prove Lemma \ref{lem:wi:ep}. 
\begin{lemma}\label{lem:evf:ep}
	We have the following convergence
\begin{align*} 
	 \inttauO{G_\ep \;[p_{{\mu}}(\vr_{\ep})]_{\ep }  \psi} \rightarrow  \inttauO{G \; \Ov{p_\mu(\vr)} \psi} \text{ as } \ep \rightarrow 0,
\end{align*}
where, $\psi \in C^1 ([0,\tau]; C^1(\T^3))$ and $\displaystyle G= \Div \vu - \Ov{\pi_{\mu}(\vr)}$.
\end{lemma}
The proof of Lemma \ref{lem:evf:ep} is performed in Appendix, in Section \ref{pf:lem:evf:ep}.  

\smallskip 

\textbf{Proof of Lemma \ref{lem:wi:ep}:} Since $\displaystyle G= \Div \vu - \Ov{\pi_{\mu}(\vr)}$, we have 
\begin{align}\label{ep:pdiv}
	&\inttauO{p_{\mu}(\vr) \Div \vu \; \psi }
	= 	\inttauO{G p_{\mu}(\vr)   \psi } + 	\inttauO{\Ov{\pi_{\mu}(\vr)} p_{\mu}(\vr)\psi}.  
\end{align}
Next, again from $\displaystyle G_{\ep}=  \Div \vu_{\ep} - [\pi_{\mu}(\vr_{\ep} )]_{\ep}$, we have
\begin{align*}
	\inttauO{p_{\mu}(\vr_{\ep} ) [\Div \vu_{\ep} ]_{\ep}  \psi }&= \inttauO{	[p_{{\mu}}(\vr_{\ep}) \psi]_{\ep } \Div \vu_{\ep}}\\
	&= \inttauO{	[p_{{\mu}}(\vr_{\ep})]_{\ep }  \Div \vu_{\ep} \psi } \\
 &\quad + \inttauO{\lr{	[p_{{\mu}}(\vr_{\ep}) \psi]_{\ep } - [p_{{\mu}}(\vr_{\ep})]_{\ep } \psi}\Div \vu_{\ep}} \\
	&=\inttauO{\lr{[p_{{\mu}}(\vr_{\ep})]_{\ep } }^2 \psi}+ \inttauO{G_\ep \;[p_{{\mu}}(\vr_{\ep})]_{\ep }  \psi} \\
	&\quad + {\inttauO{ [[\mathsf{q}(\vr_\ep)]_{\ep}]_{\ep} \;p_{{\mu}}(\vr_{\ep}) \psi} }\quad \\
	&\quad +\inttauO{\lr{	[p_{{\mu}}(\vr_{\ep}) \psi]_{\ep } - [p_{{\mu}}(\vr_{\ep})]_{\ep } \psi}\Div \vu_{\ep}}.
\end{align*}
We invoke Lemma \ref{lem:evf:ep}, Friedrich's Commutator Lemma and perform the limit passage to deduce
\begin{align}\label{ep:dpdiv}
	\inttauO{\Ov{\Ov{p_{\mu}(\vr) \Div \vu}} \;\psi } =	 \inttauO{ \mathcal{P} \; \psi }+ \inttauO{G{\Ov{\; p_{\mu}(\vr)}}  \;  \psi} + \inttauO{ \Ov{\Ov{\mathsf{q}(\vr)p_\mu(\vr)}} \psi},
\end{align}
where 
\begin{align*}
	[[\mathsf{q}(\vr_\ep)]_{\ep}]_{\ep} \;p_{{\mu}}(\vr_{\ep}) \rightharpoonup \Ov{\Ov{\mathsf{q}(\vr)p_\mu(\vr)}} \text{ in } L^{p}_{t,x} \text{ for some } p>1. 
\end{align*}
	From \eqref{ep:dpdiv} and \eqref{ep:pdiv}, we complete the proof of Lemma \ref{lem:wi:ep}. $ \square $ \par

\bigskip

Now, using identity \eqref{id:pdiv}, we reformulate equation \eqref{wfbarp-p} as
\begin{align}\label{id:phi}
	\begin{split}
		&\intO{	 \lr{\Ov{P_{\mu}(\vr)} -P_{\mu}(\vr)} (\tau) \psi(\tau) } - \inttauO{	\lr{\Ov{P_{\mu}(\vr)} -P_{\mu}(\vr)} \lr{\pt \psi + \vu \cdot \nabla \psi }} \\ 
		&+ \inttauO{G \lr{\Ov{p_{\mu}(\vr)} -p_{\mu}(\vr)  } \psi } + 	\inttauO{ \lr{\mathcal{P}- \Ov{p_{\mu}(\vr)} p_{\mu}(\vr)  }\psi  } \\
		&+ {\inttauO{ \lr{\Ov{\Ov{\mathsf{q}(\vr) p_\mu(\vr)}}  -\Ov{\mathsf{q}(\vr)} p_\mu(\vr)} \psi } }+  \inttauO{\delta  \lr{\Ov{ \vr^{m}P_{\mu}^\prime(\vr)} {-} \Ov{\vr^m}P_{\mu}^\prime(\vr)} \psi }	=0 
	\end{split}
\end{align}
for any $\psi \in C^1_{t,x}$.
In this regard, we analyse the terms of \eqref{id:phi}: 
\begin{itemize}
	\item \textbf{Sign of the term  $ \lr{\Ov{ \vr^{m}P_{\mu}^\prime(\vr)} {-} \Ov{\vr^m}P_{\mu}^\prime(\vr)} $}: \\ 
	Since both $\vr\mapsto \vr^m$ and $\vr \mapsto P'_\mu(\vr)$ are non-decreasing and $P_\mu'$ is convex, using first \cite[Theorem 11.26]{FN} and then \cite[Theorem 11.27]{FN} we have  $$ \lr{\Ov{ \vr^{m}P_{\mu}^\prime(\vr)} {-} \Ov{\vr^m}P_{\mu}^\prime(\vr)} \geq 0 .$$ 
	\item \textbf{Sign of the term $\lr{\mathcal{P}- \Ov{p_\mu(\vr)} p_\mu(\vr)  } $}: \\ 
	Since $F(\lambda) = \lambda^2 $ is convex, we obtain the following inequality in the distribution sense
	$\mathcal{P} \geq \lr{  \Ov{p_{\mu}(\vr) }   }^2.$
	Now convexity of $p_{\mu} $ gives us 
	$\Ov{p_{\mu}(\vr)}\geq p_{\mu}(\vr).$
	Hence, we have 
	\begin{align}
		\mathcal{P}\geq  \lr{  \Ov{p_{\mu}(\vr) }   }^2\geq \Ov{p_{\mu}(\vr)}\; p_{\mu}(\vr). \nonumber
	\end{align} 
 
 
	\item \textbf{Estimate for $ \lr{\Ov{\Ov{\mathsf{q}(\vr) p_\mu(\vr)}}  -\Ov{\mathsf{q}(\vr)} p_\mu(\vr)}$}: \\
	Here, we first observe that $\displaystyle \Ov{\mathsf{q}(\vr)} = \Ov{\Ov{\mathsf{q}(\vr)}}$. Next, we have 
	\begin{align}
	&\inttauO{\Big(\Ov{\Ov{\mathsf{q}(\vr) p_\mu(\vr)}} - \Ov{\mathsf{q}(\vr)} p_\mu(\vr)\Big) \psi }  \nonumber \\
  = &\lim_{\ep \to 0}\inttauO{\Big(\int_{\T^3} \mathsf{q}(\vr_\ep)(t,y) \kappa_\ep(x-y) \dy\Big) [p_\mu(\vr_\ep) -p_\mu(\vr)] \psi} ,  \nonumber
	\end{align}
	where $\kappa_{\ep}$ is the regularizing kernel in space variables. Since, $\mathsf{q}$ is bounded and compactly supported, we have 
	\[ 	\sup_{t\in (0,T)} \sup_{x\in \T^3} \Big|\int_{\T^d} \mathsf{q}(\vr_\ep(t,y)) \kappa_\ep(x-y) \dy\Big| \leq C.\]
	From \eqref{P-p-mu} we obtain
	\begin{align*}
		| \Ov{p_\mu (\vr)} - p_\mu (\vr) | \leq \lambda_{\mathsf{q}} \lr{ \Ov{P_\mu (\vr)} - P_\mu (\vr) } .
	\end{align*}
	Hence, we have
	\begin{align}\label{id:q-est}
		\left\vert {\inttauO{ \lr{\Ov{\Ov{\mathsf{q}(\vr) p_\mu(\vr)}}  -\Ov{\mathsf{q}(\vr)} p_\mu(\vr)} \psi } } \right\vert
		\leq C_{\mathsf{q}}\lambda_{\mathsf{q}} \inttauO{ \lr{\Ov{P_\mu (\vr)} - P_\mu (\vr)} \psi  }.
	\end{align}
	
	\item From \eqref{P-p-mu} we have 
	\begin{align*}
		\left\vert \inttauO{G \lr{\Ov{p_{\mu}(\vr)} -p_{\mu}(\vr)  } \psi } \right\vert \leq \inttauO{ \lambda_{\mathsf{q}} |G| \lr{\Ov{P_{\mu}(\vr)} -P_{\mu}(\vr)  } \psi }.
	\end{align*}
	
\end{itemize}
Hence, combining all the information we gathered in the previous discussion,  for $\psi\geq 0$, we obtain 
\begin{align}\label{id:phi2}
	\begin{split}
		&\intO{	 \lr{\Ov{P_{\mu}(\vr)} -P_{\mu}(\vr)} (\tau) \psi(\tau) } - \inttauO{	\lr{\Ov{P_{\mu}(\vr)} -P_{\mu}(\vr)} \lr{\pt \psi + \vu \cdot \nabla \psi }} \\ 
		&+ \inttauO{ \lr{-\lambda_{\mathsf{q}}|G|- C_{\mathsf{q}}\lambda_{\mathsf{q}}} \lr{\Ov{P_{\mu}(\vr)} -P_{\mu}(\vr)  } \psi } + 	\inttauO{ \lr{\mathcal{P}- \Ov{p_{\mu}(\vr)} p_{\mu}(\vr)  }\psi  } \\
		&+  \inttauO{\delta  \lr{\Ov{ \vr^{m}P_{\mu}^\prime(\vr)}  - \Ov{\vr^m}P_{\mu}^\prime(\vr)} \psi }	\leq 0. 
	\end{split}
\end{align} 
\begin{itemize}

	\item \textbf{Use of weight $w$: }
	Formally, substituting $\psi = w$ in \eqref{id:phi} and using \eqref{id:phi2} and \eqref{weight-eqn} yields
	\begin{align}
		\begin{split}
			&\intO{	 \lr{\Ov{P_{\mu}(\vr)} -P_{\mu}(\vr)} (\tau) w(\tau) } + \inttauO{ \lr{\Lambda-\lambda_{\mathsf{q}}|G|- C_{\mathsf{q}}\lambda_{\mathsf{q}}}\lr{\Ov{P_{\mu}(\vr)} -P_{\mu}(\vr)  } w }  \nonumber \\
			&+ 	\inttauO{ \lr{ \mathcal{P}- \Ov{p_{\mu}(\vr)} p_{\mu}(\vr)  } w  } +  \inttauO{\delta  \lr{\Ov{ \vr^{m}P_{\mu}^\prime(\vr)}  - \Ov{\vr^m}P_{\mu}^\prime(\vr)} w }	\leq 0. 
		\end{split}
	\end{align}
	We take $\Lambda \geq \lambda_{\mathsf{q}}|G|+ C_{\mathsf{q}}\lambda_{\mathsf{q}}$ and using the sign of terms along with the non-negativity of $w$ we obtain
	\begin{align}\label{58}
		\intO{	 \lr{\Ov{P_{\mu}(\vr)} -P_{\mu}(\vr)} (\tau) w(\tau) }  \leq 0 \text{ for a.e. } \tau \in (0,T).
	\end{align}
	\item Since $w$ is only in $L^\infty((0,T)\times \T^3)$, we substitute 
	$ \psi = [w]_{\delta}  $, a space-time regularization of $w$, and use commutator estimate to conclude relation (\ref{58}).

\end{itemize}


\subsubsection*{Step 2: Removal of the weight}

We intend to remove the weight from the identity
\begin{equation}\label{id-weight}
	\intO{\lr{\Ov{P_{\mu}(\vr)} -P_{\mu}(\vr)}(t) w(t) } =0
	\mbox{ \ \ for \ \ a.e. \ } t \in [0,T].
\end{equation}
From the property of the weight \eqref{boundlogw} we have that 
$\displaystyle \frac{d}{dt} \intO{ \varrho |\ln w| } \in L^1(0,T).$
We apply this information to control the weak limit $\Ov{{P_{\mu}}(\vr)}$. We proceed in the following way by considering 
\begin{equation}
	{P_{\mu}}(s)={P_{\mu}}_M(s) + b_M(s), \nonumber
\end{equation}
where 
\begin{equation*}
	{P_{\mu}}_M(s)=\left\{
	\begin{array}{cc}
		{P_{\mu}}(s) & {P_{\mu}}(s) \leq M \\
		M  &   {P_{\mu}}(s) >M
	\end{array}
	\right.
	\mbox{ \ and \ }
	b_M(s) = {P_{\mu}}(s)-{P_{\mu}}_M(s).
\end{equation*}
By \eqref{bound-n2} we know that ${P_{\mu}}(\vr_\ep)$ is uniformly bounded in $ L^{\frac{2(m-1)}{\Gamma}}((0,T)\times \T^3)$ with $ \frac{2(m-1)}{\Gamma}\geq 2$. Next, from the observation
\begin{equation*}
	{\rm supp} \; b_M(\vr_\ep) \subset \{ (t,x) \in(0,T) \times \T^3 : {P_{\mu}}(\vr_\ep) >M\},
\end{equation*}
we use the Czebyshev inequality and obtain
\begin{equation}
	|{\rm supp} \; b_M(\vr_\ep)| \leq M^{-2} \|{P_{\mu}}(\vr_\ep)\|_{L^2_{t,x}}^2 \leq C M^{-2}. \nonumber
\end{equation}
Eventually, it yields
\begin{equation*}
	\| b_M(\vr_\ep)\|_{L^{1}_{t,x}} \leq |\{ (t,x) \in  (0,T)\times \T^d: {P_{\mu}}(\vr_\ep) >M\}|^{1/2}
	\|b_M(\vr_\ep)\|_{L^2_{t,x}} \leq C M^{-1}.
\end{equation*}
Hence, we deduce
\begin{equation}\label{bm}
	\|\Ov{b_M(\vr)}\|_{L^{1}_{t,x}} \leq C M^{-1}. 
\end{equation}

Next, for the term ${P_{\mu}}_M(\vr_\ep)$ we  first observe that for some increasing function $\Phi(\cdot)$ the following holds:
\begin{equation*}
	{P_{\mu}}_M(\vr_\ep) \leq \Phi(M) \vr_\ep, \mbox{ \ since \ } {P_{\mu}}_M(0)=0.
\end{equation*}
Then for any measurable set $\omega \subset  (0,T) \times \T^d$ we have
\begin{equation*}
	\int_\omega {P_{\mu}}_M(\vr_\ep) \; \dx \ \dt \leq \Phi(M) \int_\omega \vr_\ep \ \dx\ \dt 
\end{equation*}
and eventually after performing the limit $\ep\rightarrow 0$ we obtain

\begin{equation}\label{pm-rho}
	\int_\omega \Ov{{P_{\mu}}_M(\vr)}\ \dx\; \dt \leq \Phi(M) 
	\int_\omega\vr\ \dx \ \dt.
\end{equation}
On the other hand, with the help of the relation
\begin{equation*}
	\Ov{{P_{\mu}}(\vr)} = \Ov{{P_{\mu}}_M(\vr)} + \Ov{b_M(\vr)},
\end{equation*}
for $\eta >0$, we have
\begin{equation*}
	\int_{(0,T)\times \T^3}  (\Ov{{P_{\mu}}(\vr)} - {P_{\mu}}(\vr)) \ \dx \ \dt \leq 
	\frac{1}{\eta} \int_{w > \eta} (\Ov{{P_{\mu}}(\vr)} - {P_{\mu}}(\vr))w \ \dx\ \dt + \int_{w \leq \eta} \Ov{{P_{\mu}}(\vr)} \ \dx \ \dt.
\end{equation*}
The first term disappears by \eqref{id-weight} (recall that $(\overline{P_\mu(\vr)}-P_\mu(\vr))w \geq 0$ a.e.), and the second one we treat as follows
\begin{equation*}
	\int_{w \leq \eta} \Ov{{P_{\mu}}(\vr)} \ \dx\ \dt = 
	\int_{w \leq \eta} (\Ov{{P_{\mu}}_M(\vr)} + \Ov{b_M(\vr)}) \ \dx\ \dt.
\end{equation*}
Note that the second term on the right-hand side is controlled by \eqref{bm}. 
To control the first one we use \eqref{pm-rho} and \eqref{boundlogw}. We have
\begin{equation}
	\int_{w \leq \eta} \Ov{{P_{\mu}}_M(\vr) }\ \dx\ \dt \leq 
	\frac{1}{|\ln \eta|} \Phi(M) \int_{w \leq \eta} 
	\vr |\ln w| \ \dx\ \dt\leq  \frac{C(M, E_{0,\mu})}{|\ln \eta|}. \nonumber
\end{equation}
Therefore, summing up we get
\begin{equation*}
	\int_{(0,T)\times \T^3} (\Ov{P_\mu(\vr)} - P_\mu(\vr))\; \dx\; \dt \leq CM^{-1}
	+ \frac{C(M, E_{0,\mu})}{|\ln \eta|}.
\end{equation*}
So taking large $M$, the first term is small and then taking $\eta$ sufficiently small we bound the second one. Hence for any $\epsilon>0$ we have
\begin{equation*}
	\int_{(0,T)\times \T^3}  (\Ov{P_\mu(\vr)} - P_\mu(\vr))\; \dx \; \dt \leq \epsilon.
\end{equation*}
Since the integrand is non-negative, we get the desired identity
\begin{equation*}
	\Ov{P_\mu(\vr)} = P_\mu(\vr) \mbox{ \ a.e. } (t,x) \in {(0,T)\times \T^3} .
\end{equation*}
Therefore, we conclude the first part of proof of Theorem \ref{th:mu} concerning the limit $\ep \to 0$.


\subsection{Limit passage $\delta >0$}
\subsubsection{Uniform estimates}
\begin{itemize}[leftmargin=*]
	
	\item[] We denote the sequence $(\vu_{\delta,\mu},\vr_{\delta,\mu})$ shortly by $(\vu_{\delta},\vr_{\delta})$ and then pass with $\delta \to 0$.  First, we quickly recall the estimates that we have at hand: 
	
	\item  \textbf{Energy estimate:} 
	From the momentum equation and equation for pressure, we get
	\begin{align}\label{ee-1-del}
		\frac{1}{2} \frac{d}{dt}\int_{\T^3}	\vert \vu_\delta \vert^2 \;\dx +   \frac{d}{dt}	\int_{\T^3} \Pi_{\mu}(\vr_\delta) \;\dx + 	\int_{\T^3} \vert \nabla \vu_\delta \vert^2 \; \dx + \int_{\T^3} \delta   \vr_{\delta}^{m}\Pi_{\mu}^\prime(\vr_\delta) \;\dx \leq 0.
	\end{align}
Since $\Pi_\mu(\vr_\delta) \sim P_\mu(\vr_\delta)$ and $\Pi'_\mu(\vr_\delta) \geq C(\vr_\delta^{\Gamma-1} -1)$ for $\vr_\delta \geq 1$, we have
	\begin{align}\label{bound-n1-del}
		\Vert  \vu_\delta  \Vert_{L^2_t W^{1,2}_x}+ \Vert \vu_\delta  \Vert_{L^\infty_t L^{2}_x}+ \Vert   P_{\mu}(\vr_\delta)  \Vert_{L^\infty_t L^{1}_x} + \delta \|\vr_\delta^{\Gamma-1+m}\|_{L^1_{t,x}} \leq C.
	\end{align}
	
  \item \textbf{Additional pressure estimate}: We obtain an additional pressure estimate (the proof is in Appendix in Section \ref{bog})  which reads as
 \begin{align}\label{bog-del}
		\Vert \vr_{\delta} \Vert_{L^{\Gamma+\alpha}_{t,x}} \leq C. 
	\end{align}
	with  $\alpha =\frac {13}{20} \Gamma -\frac 1{20}$ corresponding to our choice $m=\frac 52\Gamma+\frac 32$. Further, with $\Gamma > \frac{21}{13}$ we also get $\alpha >1$.

\end{itemize}
\subsubsection{Limit passage}
\begin{itemize}
	\item After limit passage the \textbf{continuity equation} 
	\begin{align*}
		&\pt \vr + \Div(\vr \vu)=0
	\end{align*}
	holds in weak sense. 
	It follows from the fact $\displaystyle \Vert \delta \vr_\delta^m \Vert_{L^{1}_{t,x}} \leq C \delta^{\frac{\Gamma-1}{\Gamma+m-1}}  \Vert \delta^{\frac{1}{\Gamma+m-1}} \vr_\delta \Vert_{L^{\Gamma+m-1}_{t,x}}^{m}$. Note that due to the square integrability of the density we also have the renormalized version of the continuity equation.
	\item \textbf{Momentum equation} after limit passage in the weak sense reads as
	\begin{align*}
		\intTO{ \vu \cdot \partial_{t} \pmb{\varphi} } + 	\intTO{  \Ov{ \pi_\mu(\vr)} \Div \pmb{\varphi} } -\intTO{  \nabla \vu : \nabla \pmb{\varphi} } =- \intO{\vu^0_\mu \cdot \pmb{\varphi}(0,\cdot)}
	\end{align*}
	\text{ for all } $\pmb{\varphi} \in C_c^1([0,T)\times \T^3) $. 
\end{itemize}
\subsubsection{Strong convergence of density}
The identification $\Ov{\pi(\vr)} = \pi(\vr) $ follows from the strong convergence of density. This holds provided we show 
\begin{align*}
	\Ov{\vr \log \vr} = \vr \log \vr \qquad \text{a.e.} 
\end{align*}
Using the renormalized form of the continuity equation for the limit problem as well as passing to the limit $\delta \to 0$ in the renormalized form on the $\delta>0$ level with the function $b(z) = z \log z$  we obtain 
\begin{align}
	\begin{split}
		\pt \lr{\Ov{\vr \log \vr} -\vr \log \vr} + \Div\lr{ \lr{\Ov{\vr \log \vr} -\vr \log \vr}   \vu }+  \lr{\Ov{\vr \Div \vu} - \vr \Div \vu } =0.\nonumber
	\end{split}
\end{align}
Note that $\displaystyle \lim_{\delta \to 0} \int_0^T\int_{\T^3} \delta \vr_\delta b'(\vr_\delta) \, \dx \,\dt = 0$.
Using the effective viscous flux equation we end up with
\begin{align}\label{ep1}
	\begin{split}
		\pt \lr{\Ov{\vr \log \vr} -\vr \log \vr} + \Div\lr{ \lr{\Ov{\vr \log \vr} -\vr \log \vr}   \vu }+  \lr{\Ov{\vr G} - \vr  G} +\lr{\Ov{\vr \pi_\mu(\vr)} - \vr \Ov{\pi_\mu(\vr)} }   =0.
	\end{split}
\end{align}
Recall that $\displaystyle G_\delta= \Div \vu_\delta - \pi_\mu(\vr_\delta) $, $\displaystyle G= \Div \vu -\Ov{ \pi_{\mu}(\vr)}$ and we have $ \vr_\delta G_\delta \rightarrow \Ov{\vr G} $ weakly in $L^{1+\lambda}$ for some $\lambda >0$.

We show below that
 \[ \Ov{\vr G} = \vr  G. \] 
First, from the uniform bound $ 	\Vert G_\delta \Vert_{L^{1{+}\frac{\alpha}{\Gamma}}_{t,x}} \leq C$ we observe that 
	$G_\delta \rightharpoonup G \text{ in } L^{1{+}\frac{\alpha}{\Gamma}}_{t,x}$.
Note that
\begin{align*}
	\inttauO{G_{\delta} \; \vr_{\delta} \psi} = -\inttauO{ \pt  (-\Delta)^{-1} \Div \vu_{\delta} \vr_\delta \psi}.
\end{align*}

Now, we introduce the space regularization of $\vu_\delta$ and denote it by $[\vu_\delta]_{\eta}$. We introduce the notation \[ {\vv_{\delta}}_\eta = \vu_{\delta}- [\vu_{\delta}]_{\eta}  \mbox{ \ \ and  \ \ } \|{\vv_{\delta}}_\eta\|_{L^{r}_{t,x}} = o(1), \quad \eta \to 0 \text{ for } r<\frac{10}{3} .\] 
Above, the  notation $o(1)$ for $ \eta \to 0$ means that the quantity goes to zero as $\eta$ goes to zero.
Therefore, we write 
\begin{align*}
	& \inttauO{  \partial_{t} {(-\Delta)}^{-1}\Div \vu_{\delta}  \;   \vr_\delta \; \psi } \\
	&= \inttauO{  \partial_{t} {(-\Delta)}^{-1}\Div [\vu_{\delta}]_{\eta}  \;   \vr_\delta \; \psi } + \inttauO{ \partial_{t} {(-\Delta)}^{-1}\Div {\vv_\delta}_\eta  \;  \vr_\delta  \; \psi }=A_1+A_2.
\end{align*}
The strong convergence of $\vu_\delta$ in $L^2_{t,x}$ and weak convergence of $\vr_\delta $ in $L^{\Gamma+\alpha}_{t,x}$ yield 
\begin{align*}
	A_1 \rightarrow \inttauO{  \partial_{t} {(-\Delta)}^{-1}\Div [\vu]_{\eta}  \;   \vr \; \psi } .
\end{align*}
On the other hand, from the uniform bound of $G_\delta$ and 
$\pt  (-\Delta)^{-1}\Div [\vu]_\eta = -[G]_\eta$,
we have 
\begin{align*}
	\inttauO{ 	\pt  (-\Delta)^{-1} \Div [\vu]_\eta \vr \;\psi} = -\inttauO{[G]_\eta \vr \; \psi}= -\inttauO{G \vr \; \psi} + o(1), \quad \eta \to 0.
\end{align*} 
For the term $A_2$ we write the equation for $( \vr_\delta \psi) $ as 
\begin{align*}
	\partial_{t}( \vr_\delta \psi) + \Div\left( \vr_\delta \vu_\delta \psi \right) 
	= -\delta \vr_\delta^{m} \psi+  \vr_\delta(\partial_t \psi + \vu_\delta \cdot \nabla  \psi) .
\end{align*}
This gives 
\begin{align*}
	\Vert \partial_{t} (  \vr_\delta \psi) \Vert_{L^{\frac 76}_t W^{-1,\frac 76}_x} \leq C(\Vert \psi \Vert_{C^1_{t,x}} ) ,
\end{align*}
provided $\Gamma$ is sufficiently large ($\Gamma \geq \frac{13}{7}$ is sufficient, as we need $\frac 76 m \leq \Gamma +m-1$).

We also have
$\Vert \vr_\delta \Vert_{L^{\Gamma+\alpha}_{t,x} }\leq C$, see Appendix \ref{bog}. Therefore we use the interpolation inequality (cf. \cite[Lemma 15]{CMZ}, we also recall it in Lemma \ref{interpolation-lemma}) with ${\vv_\delta}_\eta \in L^{\frac{10}{3}-\lambda}$, $\lambda >0$, small, and $\partial_t {\vv_\delta}_\eta\in L^{1+\frac{\alpha}{\Gamma}}_{t} W^{-1,1+\frac{\alpha}{\Gamma}}_x$. From $\Gamma $ large as above and $\alpha =\frac{13}{20}\Gamma - \frac{1}{20}$, bigger than 1, we notice 
\[ \frac{\Gamma}{\Gamma+\alpha} +\frac{1}{\Gamma +\alpha} <1 \text{ \ \ and \ \ }\frac{1}{\frac{10}{3}-\lambda} + \frac 67 >1\] and these yield
\begin{align*}
	|A_2| \leq \Vert \vr_\delta \Vert_{L^{\Gamma +\alpha}_{t,x}}^a 	\Vert \partial_{t} (  \vr_\delta \psi) \Vert_{L^{\frac 76}_t W^{-1,\frac 76}_x}^{1-a}   \Vert {\vv_\delta}_\eta \Vert_{ L^{\frac{10}{3}-\lambda}}^{1-a} \Vert \partial_t {\vv_\delta}_\eta \Vert_{ L^{\frac{\Gamma+\alpha}{\Gamma}}_{t} W^{-1,\frac{\Gamma+\alpha}{\Gamma}}_x}^a 
	\leq o(1) \mbox{ \ for $\eta \to 0$,}
\end{align*}
 $\lambda$ positive, small.
Hence, we conclude 
\begin{align*}
	\inttauO{G_{\delta}  \vr_{\delta} \psi} \rightarrow 	\inttauO{G \; \vr \psi} .
\end{align*}
Now going back to \eqref{ep1}, we have
\begin{align}\label{del-final}
	\frac{d}{dt} \intO{\lr{\Ov{\vr \log \vr} -\vr \log \vr}} + \inttauO{\lr{\Ov{\vr p_\mu(\vr)} - \vr \Ov{p_\mu(\vr)} } } =0 .
\end{align}
We notice that $\lr{\Ov{\vr \pi_\mu(\vr)} - \vr \Ov{\pi_\mu(\vr)} } = \lr{\Ov{\vr p_\mu(\vr)} - \vr \Ov{p_\mu(\vr)} }+ \lr{\Ov{\vr \mathsf{q}(\vr)} - \vr \Ov{\mathsf{q}(\vr)} }$. 

First, we observe $\displaystyle \lr{\Ov{\vr p_\mu(\vr)} - \vr \Ov{p_\mu(\vr)} }\geq 0$.
For the term $\lr{\Ov{\vr \mathsf{q}(\vr)} - \vr \Ov{\mathsf{q}(\vr)} }$ we write
\begin{align*}
	\left( \overline{\vr \mathsf{q}(\vr) } - \vr\; \Ov{\mathsf{q}(\vr)} \right) = \left( \overline{\vr \mathsf{q}(\vr) } - \vr\; {\mathsf{q}(\vr)} \right) + \left( \overline{\vr \mathsf{q}(\vr) } - \vr\; \Ov{\mathsf{q}(\vr)} \right). 
\end{align*}
Since $\mathsf{q}\in C_c^2[0,\infty) $, we fix $M>0$ such that  $\text{supp} \; \mathsf{q} \subset [0,M)$.   Therefore, there exists a constant $\lambda_{\mathsf{q}}>0$ that depends on $\mathsf{q}$ such that 
\begin{align*}
	\lambda_{\mathsf{q}} \vr \log \vr +\vr \mathsf{q}(\vr) , 	 \quad \lambda_{\mathsf{q}} \vr \log \vr +\mathsf{q}(\vr)  \text{ \ \ and \ \ } 	\lambda_{\mathsf{q}} \vr \log \vr - \vr \mathsf{q}(\vr)  \text{ \ \ are convex}.
\end{align*}
The above follows from the fact that 
	$\mathsf{q}\in C_c^2[0,\infty) , \sup_{(m,M)} \lr{|\mathsf{q}| +|\mathsf{q}^\prime| + |\mathsf{q}^{\prime \prime}|} <\infty .$
This helps us to deduce 
\begin{align*}
	\left\vert \int_0^\tau \int_{\T^3} \left( \overline{\vr \mathsf{q}(\vr) } - \vr\; \Ov{\mathsf{q}(\vr)} \right) \, \dx\, \dt \right\vert \leq C(\mathsf{q},\lambda_\mathsf{q})\int_0^\tau  \int_{\T^3} \left(\overline{\vr \log \vr} -\vr \log \vr \right)(\tau,\cdot )\, \dx \,\dt. 
\end{align*}
Finally, we use the Gr\"onwall's argument in \eqref{del-final} to conclude
	$\Ov{\vr \log \vr} =\vr \log \vr$.
 Thus, up to a subsequence $\vr$ convergences point-wisely. Theorem \ref{th:mu} is proved.


\section{Proof of Theorem \ref{th:main}}

The statement of the theorem outlines two cases of pressures related to the growth condition. The first one is dedicated to the case with compact perturbation of a monotone pressure, i.e.,
\begin{align*}
     &\pi(\vr)= p(\vr)+ q(\vr),\; \text{ with } q\in C_c[0,\infty), \; q(0)=0 \text{ and } \\
     & p(0)=0,\; p(\vr) > 0,\; p^\prime(\vr) >0 \text{ for } \vr>0 \text{ and } p \text{ satisfies } \eqref{pre1} \text{ with }\gamma >1.\nonumber
     \qquad {{\bf \rm(Case\;1)}}
 \end{align*} 
The second one is for the non-monotone case with certain growth hypotheses. We take a smooth non-negative function $p\in C^1([0,\infty))\cap C^2((0,\infty))$ such that $p(\vr) >0$ for $\vr >0$ and 
\begin{equation*}
	p(0)=0 \text{  \ \  and  \ \ }  a_2 \vr^\gamma -C \leq p(\vr) \leq C + a_1 \vr^\gamma \mbox{ \ \ with \ \ } \gamma \geq 6/5,\; a_1,a_2,C>0.
  \qquad {{\bf \rm(Case\;2)}}
\end{equation*}
Additionally, we require 
\begin{equation*}
	|p'(\vr)|\leq C \vr^{\gamma-1}, \qquad  |p''(\vr)|\leq C \vr^{\gamma-2} \mbox{ \ \ for \ \ } \vr > 1.
\end{equation*}


\subsection{Uniform estimates}
\begin{itemize}[leftmargin=*]
	
	\item[] The main goal here is to perform the limit passage $\mu \rightarrow0$. To achieve this, first we recall the uniform estimates w.r.t $\mu$: 
	
	\item  \textbf{Energy estimate:} 
	From the energy inequality, for a.e. $\tau \in (0,T)$ we have
	\begin{align}
		&\int_{\T^3}	\lr{\frac{1}{2} \vert \vu_\mu \vert^2 \;\dx + \Pi_\mu(\vr_\mu)\;} (\tau) \dx + 	\int_0^\tau \int_{\T^3} \vert \nabla \vu_\mu \vert^2 \; \dx \; \dt 
		&\leq 	\int_{\T^3}	\lr{\frac{1}{2} \vert \vu_{\mu,0} \vert^2 \;\dx + \Pi_\mu(\vr_{\mu,0})\;}   \dx. \nonumber
	\end{align}
	This yields
	\begin{align}
		\Vert  \vu_\mu  \Vert_{L^2_t W^{1,2}_x}+ \Vert \vu_\mu  \Vert_{L^\infty_t L^{2}_x}+ \Vert   \vr_\mu^\gamma \Vert_{L^\infty_t L^{1}_x} +  \Vert  \mu \vr_\mu^\Gamma \Vert_{L^\infty_t L^{1}_x}\leq C. \nonumber
	\end{align}
	\item \textbf{Extra pressure estimate}: We obtain an additional pressure estimate (details are in Appendix \ref{bog}, case 2) 
	\begin{align}\nonumber
		\Vert \vr^{\frac{5}{3} \gamma}_{\mu} \Vert_{L^1_{t,x}} \leq C. 
	\end{align}
	In particular, for $\gamma> 1$ we get
		$\Vert p(\vr_{\mu}) \Vert_{L^p_{t,x}} \leq C \text{ with } p >1. $

\end{itemize}
\subsection{Limit passage}
\begin{itemize}[leftmargin=*]
	\item After the limit passage, the \textbf{continuity equation} 
	\begin{align*}
		&\pt \vr + \Div(\vr \vu)=0
	\end{align*}
	holds in a weak sense. 
	\item After the limit passage, the \textbf{velocity equation} in weak sense reads as
	\begin{align}
		\iintTOM{ \vu \cdot \partial_{t} \pmb{\varphi} } + 	\iintTOM{  \Ov{ \pi(\vr)} \Div \pmb{\varphi} } -\iintTOM{  \nabla \vu : \nabla \pmb{\varphi} } =- \intO{\vu_0 \cdot \pmb{\varphi}(0,\cdot)} \nonumber
	\end{align}
	\text{ for all } $\pmb{\varphi} \in C_c^1([0,T)\times \T^3) $.
	\item \textbf{Renormalized equation of continuity} holds for the limit system. The problematic situation is for $\gamma \le \frac 65$ and it can be achieved only for the pressure with compact non-monotone perturbation. The proof is nontrivial and needs the concept of oscillation defect measure introduced by Feireisl \cite{EF2001}. 
	
\end{itemize}

\subsection{Strong convergence of density}
\subsubsection{Case 1: Pressure satisfying (\ref{pr-cmpt})}
The identification $\Ov{\pi(\vr)} = \pi(\vr)$ follows from the strong convergence of density, specifically from $\displaystyle \Ov{\vr \log \vr} = \vr \log \vr$. The key idea is based on the work of Feireisl \cite{F2002}. Below, we provide a list of important propositions and statements used to prove the strong convergence of density.
\begin{itemize}[leftmargin=*]
	\item \textbf{Effective viscous flux identity:} 
	\begin{prop}\label{Prop-EVF-b}
		For suitable $b$ (the main points are that the function cannot be too singular at zero and cannot grow too fast at infinity) we have that the following effective viscous flux identity holds in the weak sense:
		\begin{align*}
			\overline{\pi(\vr) b(\vr) }- \overline{\pi(\vr)}\;\overline{b(\vr)} =\overline{b(\vr)\dv \vu} - \overline{b(\vr)} \dv \vu .
		\end{align*}
  \end{prop} 
  
  The proof follows the same lines as in the case of compressible Navier--Stokes equations. However, since the momentum equation does not include the convective term and we have directly the strong convergence of the velocity, it is more direct as it does not require any compensated compactness tools.
	
	\item \textbf{Oscillation defect measure}: We briefly recall this concept. First, we consider a smooth cut-off function $ T $: 
	\begin{align}\label{Tt}
		T(x) = \begin{cases}
			&x,\; 0\leq x\leq 1\\
			& \text{concave smooth},  1\leq x\leq 3\\
			& 2,\; x \geq 3.\\
		\end{cases}
	\end{align}
	Then we define ``Truncation''  functions: 
		$T_k(x) = k \;T\left( \frac{x}{k}\right).$
	For $\alpha>0$, we define \emph{Oscillation defect measure} as
	\begin{align*}
		\textbf{osc}_\alpha [\vr_\mu \rightarrow \vr ]  := \sup_{k\geq 1} \left(\limsup_{\ep\to0} \int_0^T \int_{\T^3} \vert T_k(\vr_\mu)- T_k(\vr)  \vert^\alpha \,\dx \,\dt  \right).
	\end{align*}
	
	\item \textbf{Boundedness of oscillation defect measure: }
	\begin{prop}\label{prop:OscD}
		Suppose, the \emph{renormalized continuity equation} and the \emph{effective viscous flux} identity is true for approximating function $\{\vr_\mu,\vu_\mu\}$, then for $ \gamma >1 $, we have 
		\begin{align*}
			\textbf{osc}_{\gamma+1} [\vr_\mu \rightarrow \vr ] <\infty . 
		\end{align*}
	\end{prop}
	The proof of the proposition follows directly from the Feireisl \cite{F2002}. For the sake of completeness it can be found in the Appendix, Section \ref{a5}.

	\item \textbf{Renormalized continuity equation: }
	\begin{prop}\label{prop:REC}
		Suppose, a sequence $ \{\vr_\mu,\vu_\mu\}$ satisfies renormalized continuity equation and they have the following convergence: \begin{align*}
			&	\vr_\mu \rightarrow \vr \text{ weak-(*)ly  in } L^\infty(0,T;L^\gamma({\T^3})),\\
			&	\vu_\mu \rightarrow \vu  \text{ weakly in } L^2((0,T)\times {\T^3}),\\
			&	\nabla \vu_\mu \rightarrow \nabla \vu  \text{ weakly in } L^2((0,T)\times {\T^3}).
		\end{align*}
		Then, if we assume they have finite oscillation defect measure, i.e.,
		\begin{align*}
			\textbf{osc}_{r} [\vr_\mu \rightarrow \vr ] <\infty  
		\end{align*}
for some $r>2$, then  $ (\vr, \vu) $ is also a solution of the renormalized continuity equation.
	\end{prop}
	The proof of the Proposition \eqref{prop:REC} is exactly the same as in the case of the compressible Navier--Stokes equations and is available e.g. in Feireisl, Karper and Pokorn\'y \cite{FKP}. For the sake of completeness, it is presented in Appendix, Section \ref{a6}.
 
	\item As a consequence of Propositions \ref{prop:OscD} and \ref{prop:REC}, we state the desired strong convergence result: 
	\begin{prop}\label{Prop-stronconv}
		Suppose, a sequence $ \{\vr_\mu,\vu_\mu\}$ satisfies renormalized continuity equation along with energy estimate, pressure estimate, and effective viscous flux identity. Then, we claim \begin{align*}
			\overline{\vr \log \vr}=\vr \log \vr.
		\end{align*}
		As a consequence of it we have 
		$ \vr_\mu \rightarrow \vr \text{ in } L^1((0,T)\times {\T^3}) \text{ as } \mu \rightarrow 0.$
	\end{prop}
\end{itemize}
\textit{Proof:} Recall that our pressure has the form 
$\pi \approx \vr^\gamma + \tilde{p}(\vr)+ \tilde q(\vr),$
where $q$ has a compact support.

We consider 
\begin{align*}
	L_k(\vr)=  \vr \int_0^\vr \frac{T_k(\xi)}{\xi^2} \,\text{d}\xi.
\end{align*}
It has the property\begin{align*}
	&	L_k^{\prime \prime}(\vr) \geq 0, \text{ for } \vr \geq 0, \qquad
		\vr L_k^\prime (\vr) - L_k(\vr)= T_k(\vr)\text{ for } \vr \geq 0.
\end{align*} 
Moreover, we notice that  $L_k(\vr) \approx \vr \log \vr \text{ as } k\rightarrow \infty.$
Seeing that 
\begin{align*}
	\frac{\text{d}}{\dt} \int_{\T^3} L_k(\vr_\mu) \,\dx + \int_{\T^3} T_k(\vr_\mu)\dv \vu_\mu \,\dx =0 \text{ for a.e. } t\in (0,T),
\end{align*}
we obtain
\begin{align*}
	\frac{\text{d}}{\dt} \int_{\T^3} \overline{L_k(\vr)} \,\dx + \int_{\T^3} \overline{T_k(\vr)\dv \vu} \,\dx =0 \text{ for a.e. } t\in (0,T).
\end{align*}
Now, from Proposition \ref{prop:REC} we conclude that $ (\vr,\vu) $ satisfies renormalized equation of continuity and it gives
\begin{align*}
	\int_{\T^3} \left(\overline{L_k(\vr)} - L_k(\vr) \right)(\tau,\cdot ) \,\dx + \int_0^\tau \int_{\T^3}\left( \overline{T_k(\vr)\dv \vu} - {T_k(\vr)} \dv \vu \right)\,\dx\,\dt =0 \text{ for a.e. } \tau \in (0,T).
\end{align*}
Next, we easily see that
\begin{align*}
	&\int_0^\tau \int_{\T^3}\left( \overline{T_k(\vr)\dv \vu} - {T_k(\vr)} \dv \vu \right)\,\dx\,\dt \\
	&=\int_0^\tau \int_{\T^3}\left( \overline{T_k(\vr)\dv \vu} - \Ov{T_k(\vr)} \dv \vu \right)\,\dx\,\dt +\int_0^\tau \int_{\T^3}\left( \overline{T_k(\vr) } - {T_k(\vr)} \right)\dv \vu \,\dx\,\dt .
\end{align*}
Using the effective viscous flux identity
\begin{align*}
	\overline{\pi(\vr) T_k(\vr) }- \overline{\pi(\vr)}\;\overline{T_k(\vr)} =\overline{T_k(\vr)\dv \vu} - \overline{T_k(\vr)} \dv \vu ,
\end{align*}
 we obtain 
\begin{align}\label{osd-1}
	&\int_{\T^3} \left(\overline{L_k(\vr)} - L_k(\vr) \right)(\tau,\cdot ) \,\dx + \int_0^\tau \int_{\T^3} \left( 	\overline{\pi(\vr) T_k(\vr) }- \overline{\pi(\vr)}\;\overline{T_k(\vr)} \right)\, \dx\,\dt \nonumber \\
	&=-\int_0^\tau \int_{\T^3}\left( \overline{T_k(\vr) } - {T_k(\vr)} \right)\dv \vu \,\dx\,\dt. 
\end{align}

First we notice that
\begin{align*}
	&\int_0^\tau \int_{\T^3}\left( \overline{T_k(\vr) } - {T_k(\vr)} \right)\dv \vu \,\dx\,\dt \\
	&\leq \Vert  {\dv \vu} \Vert_{L^2((0,T)\times {\T^3}) }  \Vert  \overline{T_k(\vr)} - {T_k(\vr)}  \Vert_{L^2((0,T)\times {\T^3}) }\\
	&\leq \Vert  {\dv \vu} \Vert_{L^2((0,T)\times {\T^3}) }  \Vert  \overline{T_k(\vr)} - {T_k(\vr)}  \Vert_{L^1((0,T)\times {\T^3}) }^\lambda \Vert \overline{T_k(\vr)} - {T_k(\vr)}  \Vert_{L^{\gamma+1}((0,T)\times {\T^3}) }^{1-\lambda}.
\end{align*}
In this case we observe (\ref{Tt}) that
\[ \Vert  \overline{T_k(\vr)} - {T_k(\vr)}  \Vert_{L^1((0,T)\times {\T^3}) } \rightarrow 0 \text{ as } k\rightarrow \infty \]
and Proposition \ref{prop:OscD}
\[\Vert \overline{T_k(\vr)} - {T_k(\vr)}  \Vert_{L^{\gamma+1}((0,T)\times {\T^3}) } <\infty .\]

\begin{itemize}
	\item \textbf{Monotone case: $q(\vr) \equiv 0$}: Using \cite[Theorem 11.26]{FN}, we have 
	\begin{align*}
		\int_0^\tau \int_{\T^3}\left( \overline{T_k(\vr) \pi(\vr) } - \Ov{T_k(\vr)}\; \Ov{\pi(\vr)} \right)\,\dx\,\dt \geq 0.
	\end{align*}
 
	Hence passing limit $k\rightarrow \infty  $ in \eqref{osd-1}, we  have 
	\begin{align*}
		 \int_{\T^3} \left(\overline{\vr \log \vr} -\vr \log \vr \right)(\tau,\cdot ) \,\dx \leq 0 \text{ for a.e. } \tau \in (0,T).
	\end{align*}
        \item \textbf{Non-monotone case with compact perturbation ($\pi \approx \vr^\gamma + \tilde{p}(\vr)+ \tilde q(\vr) $)}: We rewrite \eqref{osd-1} as
\begin{align*}
	&\int_{\T^3} \left(\overline{L_k(\vr)} - L_k(\vr) \right)(\tau,\cdot ) \,\dx + \int_0^\tau \int_{\T^3} \left( 	\overline{\lr{\vr^\gamma+\tilde{p}(\vr)}T_k(\vr) }- \overline{\lr{\vr^\gamma+\tilde{p}(\vr)}}\;\overline{T_k(\vr)} \right)\,\dx\,\dt \\
	&=-\int_0^\tau \int_{\T^3}\left( \overline{T_k(\vr) } - {T_k(\vr)} \right)\dv \vu \,\dx\,\dt - \int_0^\tau \int_{\T^3} \left( 	\overline{\tilde{q}(\vr) T_k(\vr) }- \overline{\tilde{q}(\vr)}\;\overline{T_k(\vr)} \right)\,\dx\,\dt \\
 &=\int_0^\tau \int_{\T^3}\left( \overline{T_k(\vr) } - {T_k(\vr)} \right)\lr{-\dv \vu+ 	\overline{\tilde{q}(\vr)} } \,\dx\,\dt - \int_0^\tau \int_{\T^3} \left( 	\overline{\tilde{q}(\vr) T_k(\vr) }- \overline{\tilde{q}(\vr)}\;{T_k(\vr)} \right)\,\dx\,\dt. 
\end{align*} 
Using the monotonicity of $\vr^\gamma+ \tilde{p}(\vr)$, we obtain sign for the second term in left-hand side. The first term in the right-hand side can be controlled in a similar way as in the monotone case, i.e.,
\begin{align*}
	&\int_0^\tau \int_{\T^3}\left( \overline{T_k(\vr) } - {T_k(\vr)} \right)\lr{-\dv \vu+ 	\overline{\tilde{q}(\vr)} } \,\dx\,\dt \\
	&\leq \Vert  \lr{-\dv \vu+ 	\overline{\tilde{q}(\vr)} } \Vert_{L^2((0,T)\times {\T^3}) }  \Vert  \overline{T_k(\vr)} - {T_k(\vr)}  \Vert_{L^1((0,T)\times {\T^3}) }^\lambda \Vert \overline{T_k(\vr)} - {T_k(\vr)}  \Vert_{L^{\gamma+1}((0,T)\times {\T^3}) }^{1-\lambda}.
\end{align*}
We again observe that
$\Vert  \overline{T_k(\vr)} - {T_k(\vr)}  \Vert_{L^1((0,T)\times {\T^3}) } \rightarrow 0$
and $\Vert \overline{T_k(\vr)} - {T_k(\vr)}  \Vert_{L^{\gamma+1}((0,T)\times {\T^3}) } <\infty .$
Therefore letting $k\rightarrow \infty $
\begin{align*}
	\int_{\T^3} \left(\overline{\vr \log \vr} -\vr \log \vr \right)(\tau,\cdot ) \,\dx \leq \left\vert \int_0^\tau \int_{\T^3} \left( \overline{\vr q(\vr) } - \vr\; \Ov{q(\vr)} \right) \,\dx\,\dt \right\vert \text{ a.e. } \tau \in (0,T).
\end{align*}
Since, $0<q\in C_c^2(0,\infty) $ we assume $\text{supp} \; q \subset [0,M]$.   Therefore, there exists a constant $\lambda(q)>0$ that depends on $q$ such that 
	\begin{align*}
		\lambda_q \vr \log \vr +\vr q(\vr) , 	\quad \lambda_q \vr \log \vr +q(\vr)  \text{ \ \ and \ \ } 	\lambda_q \vr \log \vr - \vr q(\vr)  \text{ \ \  are convex}.
	\end{align*}
 This helps us to deduce 
	\begin{align*}
		\left\vert \int_0^\tau \int_{\T^3} \left( \overline{\vr q(\vr) } - \vr\; \Ov{q(\vr)} \right) \,\dx\,\dt \right\vert \leq C(q,\lambda_q)\int_0^\tau  \int_{\T^3} \left(\overline{\vr \log \vr} -\vr \log \vr \right)(\tau,\cdot ) \,\dx \,\dt. 
	\end{align*}
Finally using Gronwall's argument we conclude
\[  \int_{\T^3} \left(\overline{\vr \log \vr} -\vr \log \vr \right)(\tau,\cdot ) \,\dx \leq 0 \text{ a.e. } \tau \in (0,T). \]
Since we already have $ \overline{\vr \log \vr}\geq \vr \log \vr $, the above inequality implies
\begin{align*}
	\overline{\vr \log \vr}=\vr \log \vr \text{ \ \   for a.e. } (t,x) \in (0,T)\times {\T^3}.
\end{align*}
This completes the proof of strong convergence.

\end{itemize}


\subsubsection{Case 2: Pressure satisfying (\ref{pre1})--(\ref{pre2}) with $\gamma \geq \frac{6}{5}$}

Here we  deal with the case on the non-monotone pressure. Our sketch of the proof will be based on the Bresch--Jabin technique introduced in \cite{BJ} and \cite{BJ2}. We need to deal with the quantity from (\ref{criterion}). In order to have sufficient differentiability of 
\[ \mathcal{K}_h(x-y)|\vr_\mu(t,x)-\vr_\mu(t,y)| ,\]
we introduce the modification of the modulus function.
Given $\sigma>0$ we define 
\eq{\nonumber
|w|^\sigma= \left\{ 
\begin{array}{lcr}
|w|-\frac{\sigma}{2} & \mbox{ for } & |w| > \sigma, \\
\frac{1}{2\sigma} w^2 & \mbox{ for } & |w|\leq \sigma,
\end{array}
\right.
}
along with
\begin{equation*}
   \sgn^\sigma:= \partial |w|^\sigma = \left\{ 
    \begin{array}{lcr}
    {\rm sgn \, }w  & \mbox{for} & |w|>\sigma \\
    \frac{w}{\sigma} & \mbox{for} & |w|< \sigma
    \end{array}
    \right. 
\mbox{ \ and \ } 
\partial^2 |w|^\sigma = \left\{ 
    \begin{array}{lcr}
    0 & \mbox{for} & |w|>\sigma \\
    \frac{1}{\sigma} & \mbox{for} & |w|\leq \sigma.
    \end{array}
    \right. 
\end{equation*}
Considering  $ \sgn^\sigma_{xy}:=\sgn^\sigma(\vr(x)-\vr(y))$, we note a key feature of this modification, which is
\begin{equation}
  |(\vr(x)-\vr(y))\sgn^\sigma_{xy} -|\vr(x)-\vr(y)|^\sigma |  \leq \sigma.\nonumber
\end{equation}
To simplify the notation we drop the dependence on parameter $\mu$ in the sequences $\{\vr_\mu,\vu_\mu\}_{\mu>0}$.  
Multiplying the difference of the continuity equation at point $x$ and $y$ by $\sgn^\sigma_{xy}$, we obtain
\eq{\label{eq:drho1}
&\pt |\vr(x)-\vr(y)|^\sigma + \dv_x (\vu(x)|\vr(x)-\vr(y)|^\sigma) +
\dv_y (\vu(y)|\vr(x)-\vr(y)|^\sigma) \\[5pt]
&+ [(\vr(x)-\vr(y))\sgn^\sigma_{xy} -|\vr(x)-\vr(y)|^\sigma](\dv_x\vu(x)+
\dv_y\vu(y))\\[5pt]
&=\dv_x \vu(x) (\vr(x)-\vr(y))\sgn^\sigma_{xy} 
 - \left[(\dv_x \vu(x)-\dv_y \vu(y))\sgn^\sigma_{xy}\right] \vr(x).
}

Here again, we need to consider the weight function $w$ defined by the solution to 
\begin{equation*}
    \partial_t w +\vu \cdot \nabla w + \Lambda w=0, \qquad w|_{t=0}=1.
\end{equation*}
In the following computation, we drop the index $\mu$ when no confusion can arise. The function $\Lambda$ is positive and satisfies $\Lambda > c_1 {\cal M}[|\Grad\vu|] + c_2 {\cal M}[\vr^\gamma] +c_3 $, where $c_1,\; c_2$ and $c_3$ are positive constants to be fixed later in the course of the proof. \par 

Also, we assume first $0\leq \vr_0 \in  L^\infty(\T^3)$ and eventually to we relax this hypothesis at the end of the proof. The first observation from the assumption $0\leq \vr_0 \in  L^\infty(\T^3)$  is 
\begin{equation}\label{rho-bdd}
\underset{t>0}{{\rm ess}\sup} \|\vr w\|_{L^\infty_x} \leq \| \vr_0\|_{L^\infty}.
\end{equation}
Using $\mathcal{K}_h(x-y) w(x)$ as a test function in \eqref{eq:drho1} and eventually using the \eqref{weight-eqn}, we get
\eqh{
\frac{d}{dt} R_h^{\sigma}(t) = A_1 + A_2 + A_3+A_4 ,
}
where 
\begin{align*}
    &R_h^{\sigma}(t) := \iintO{\mathcal{K}_h(x-y) |\vr(x)-\vr(y)|^\sigma w(x)},
\end{align*}
and 
\begin{align*}
&A_1= 
    -\iintO{ \mathcal{K}_h(x-y) [\Div_x (\vu(x) |\vr(x) -\vr(y)|^\sigma w(x))
    +\Div_y (\vu(y)|\vr(x)-\vr(y)|^\sigma w(x)] }
    \\
&A_2=     
    \iintO{ \mathcal{K}_h(x-y) (\Div_x \vu(x) - \Lambda) |\vr(x)-\vr(y)|^\sigma  w(x) }\\
&A_3=     - \iintO{ \mathcal{K}_h(x-y)(\Div_x \vu(x) - \Div_y \vu(y)) \sgn_{xy}^\sigma  \;\vr(x) w(x) } \\
&A_4=     -\iintO{ \int \mathcal{K}_h
    [(\vr(x)-\vr(y))\sgn^\sigma_{xy} -|\vr(x)-\vr(y)|^\sigma](2\dv_x\vu(x)+
\dv_y\vu(y)) w(x) }.
\end{align*}

The last term is just bounded by $\sigma$, in particular
\begin{equation}\label{a4}
    \int_0^T |A_4|  \;{\rm d}t \leq C \sigma,
\end{equation}
where the constant $C$ depends only on the initial data.
For the term $A_1$, we have 
\begin{align*}
     &-\iintO{ \mathcal{K}_h(x-y) [\Div_x (\vu(x) |\vr(x) -\vr(y)|^\sigma w(x))
    +\Div_y (\vu(y)|\vr(x)-\vr(y)|^\sigma w(x))] } \\
    &=\iintO{\nabla \mathcal{K}_h(x-y)\lr{ \vu(x)-\vu(y)} |\vr(x)-\vr(y)|^\sigma w(x)}.
\end{align*}
Therefore, using the property of the kernel $\mathcal{K}_h$
\[|z||\Grad \mathcal{K}_h(z)|\leq C \mathcal{K}_h(z),\]
we obtain 
\eq{\nonumber
|A_1|&\leq \iintO{|\nabla \mathcal{K}_h(x-y)|| \vu(x)-\vu(y)| |\vr(x)-\vr(y)|^\sigma w(x)}\\
&\leq C \iintO{\mathcal{K}_h(x-y) \frac{|\vu(x)-\vu(y)|}{|x-y|} |\vr(x)-\vr(y)|^\sigma w(x)}.
}
Here, we recall the following lemma (see \cite{BJ} and \cite[Lemma 11]{CMZ})
\begin{lemma}\label{Lagrange}
Let $f \in W^{1,p}(\T^d)$ with $p>1$, then 
\begin{equation}\label{formula:D}
    |f(x)-f(y)|\leq C|x-y|\left( D_{|x-y|} f(x) + D_{|x-y|}f(y)\right) \mbox{ \ a.e. } {x,y \in \T^d,}
\end{equation}
where for $r>0$ we denoted
\begin{equation}\nonumber
D_rf(x)=\frac{1}{r}\int_{B(0,r)} \frac{|\nabla f(x+z)|}{|z|^{d-1} } {\rm d}z.
\end{equation}
\end{lemma}
Using the above lemma, we have 
\eqh{
|A_1|\leq& C\iintO{\mathcal{K}_h(x-y) \lr{D_{|x-y|}\vu(x)+D_{|x-y|}\vu(y)} |\vr(x)-\vr(y)|^\sigma w(x)}\\
\leq& C\iintO{\mathcal{K}_h(x-y) \lr{D_{|x-y|}\vu(x)-D_{|x-y|}\vu(y)} |\vr(x)-\vr(y)|^\sigma w(x)}\\
&+C\iintO{\mathcal{K}_h(x-y) {\cal M}[|\Grad\vu(x)|] |\vr(x)-\vr(y)|^\sigma w(x)}.
}
Recall that 
\[A_2=     
    \iintO{ \mathcal{K}_h(x-y) (\Div_x \vu(x) - \Lambda) |\vr(x)-\vr(y)|^\sigma  w(x) }.\]
Since, ${\cal M}[|\Grad\vu|]$ is bounded in $L^2_{t,x}$, it yields 
\begin{align}\label{a1a2}
\begin{split}
    |A_1|+ A_2\leq& C\iintO{\mathcal{K}_h(x-y) \lr{D_{|x-y|}\vu(x)-D_{|x-y|}\vu(y)} |\vr(x)-\vr(y)|^\sigma w(x)} \\
    & + C\iintO{ \mathcal{K}_h(x-y) ( c_1 {\cal M}[|\Grad\vu|]- \Lambda) |\vr(x)-\vr(y)|^\sigma  w(x) },
    \end{split}
\end{align}
for some $c_1>0$.

\textbf{Estimate for the term $A_3$:} This is one of the most crucial part of the analysis. It follows from the nice features of effective viscous flux. We recall the effective viscous flux 
\[ G_\mu= \Div \vu_\mu- p_\mu(\vr_\mu). \] 

Using this notation (dropping the index $\mu$ for the terms that are uniformly estimated), we have
\begin{align*}
    A_3=   &  - \iintO{ \mathcal{K}_h(x-y)(\Div_x \vu(x) - \Div_y \vu(y)) \sgn_{xy}^\sigma  \;\vr(x) w(x) } \\
    =&- \iintO{ \mathcal{K}_h(x-y)(G(x) - G(y)) \sgn_{xy}^\sigma  \;\vr(x) w(x) } \\
    &-\iintO{ \mathcal{K}_h(x-y)(p_\mu(\vr_\mu(x)) - p_\mu(\vr_\mu(y))) \sgn_{xy}^\sigma  \;\vr(x) w(x) }\\
    =&- \iintO{ \mathcal{K}_h(x-y)(G(x) - G(y)) \sgn_{xy}^\sigma  \;\vr(x) w(x) } \\
    &-\iintO{ \mathcal{K}_h(x-y)(p(\vr(x)) - p(\vr(y))) \sgn_{xy}^\sigma  \;\vr(x) w(x) }\\
    &-\mu \iintO{ \mathcal{K}_h(x-y)(\vr^{\Gamma}(x) - \vr^{\Gamma}(y)) \sgn_{xy}^\sigma  \;\vr(x) w(x) }\\
    =&A_{3,1}+A_{3,2}+A_{3,3}.
\end{align*}

The estimate of the term $A_{3,1}$ in the right-hand side is a consequence of application of the the interpolation inequality, so first we look at the terms $A_{3,2}$ and $A_{3,3}$.
Using the monotonicity of $\lambda \mapsto \lambda^\Gamma$, we have the desired sign for the term $A_{3,3}$. 
Next, for the term $A_{3,2}$ we follow strategy from \cite{BJ2}, we first consider the regime $\left\{ (x,y)\in \T^6 \;\big| \;  \vr(y) \leq \vr(x) \right\}$. Eventually,  we consider two sub-regimes: 
\begin{align*}
    &1. \;\Omega_1=\left\{ (x,y)\in \T^6 \;\bigg| \;  \vr(y) \leq \frac{\vr(x)}{2}  \right\},\\
    &2. \;\Omega_2=\left\{ (x,y)\in \T^6 \;\bigg| \;  \frac{\vr(x)}{2} \leq \vr(y) \leq \vr(x)  \right\}.
\end{align*}
In $\Omega_2$, we use the Lagrange's lemma (Lemma \ref{Lagrange}) to obtain
\begin{equation*}
    (p(\vr(x)-p(\vr(y))\sgn_{xy}^\sigma \geq -C(\vr_x^{\gamma-1} +1)|\vr(x)-\vr(y)|.
\end{equation*}
On the other hand, in $\Omega_1$, at first we note that $\vr(y) \leq \vr(x)/2$ implies 
\begin{equation}\nonumber
    p(\vr(x))-p(\vr(y)) \geq - C((\vr(y))^\gamma +1) \geq -C(\vr(x))^\gamma +1)
\end{equation}
and 
\begin{equation}\nonumber
    \vr(x) \leq \vr(x) + \vr(x) - 2\vr(y)\leq 
    2(\vr(x)-\vr(y)).
\end{equation}
This yields
\begin{equation*}
    [p(\vr(x)) - p(\vr(y))] \sgn_{xy}^\sigma \vr(x)
    \geq -C(\vr(x)^\gamma +1)|\vr(x) -\vr(y)|.
\end{equation*}
Using a similar analysis we have the same bound for the case $\left\{ (x,y)\in \T^6 \;\big| \;  \vr(y) \leq \vr(x) \right\}$.  Thus, we deduce 
\begin{align*}
    A_{3,2} \leq  \iintO{ \mathcal{K}_h(x-y) ( c_2 {\cal M}[\vr^\gamma+1]) |\vr(x)-\vr(y)|^\sigma  w(x) } + C \sigma. 
\end{align*}
Therefore, the above estimate along with \eqref{a4} and \eqref{a1a2} gives
\begin{align}\label{a-a31}
    &|A_1|+ A_2 + A_{3,2}+ A_{3,3}+ A_4 \nonumber \\ 
    &\leq C\iintO{\mathcal{K}_h(x-y) \lr{D_{|x-y|}\vu(x)-D_{|x-y|}\vu(y)} |\vr(x)-\vr(y)|^\sigma w(x)} \nonumber\\
    & \quad + \iintO{ \mathcal{K}_h(x-y) ( \lr{c_1 {\cal M}[|\Grad\vu|] + c_2 {\cal M}[\vr^\gamma] +c_3} - \Lambda) |\vr(x)-\vr(y)|^\sigma  w(x) } + C \sigma,
\end{align}
for some $c_1,c_2, c_3>0$.
Therefore, for some $\Lambda > c_1 {\cal M}[|\Grad\vu|] + c_2 {\cal M}[\vr^\gamma] +c_3 $, we can ignore the term in the right-hand side.\par 
\medskip 

Hence, the remaining uphill task is to estimate the term $A_{3,1}$.
First we rewrite term as
\begin{align*}
   A_{3,1}= \iintO{ \mathcal{K}_h(x-y)\lr{\pt (\Delta_x)^{-1} \dv_x \vu(x) - \pt (\Delta_y)^{-1} \dv_y\vu(y)} W },
\end{align*}
where 
$W=W(t,x,y)= \sgn^\sigma_{xy} \, \vr(x)   w(x)$
 and the equation for effective viscous flux \eqref{if:evf1} that reads as
     $G= -\pt (-\Delta)^{-1} \dv \vu.$

Let $\kappa_\eta$ be the space regularization kernel and we denote the space regularization of $\vu$ by $[\vu]_\eta$. Moreover we denote 
\begin{align*}
   \vv_\eta = \vu -[\vu]_\eta. 
\end{align*}
Using this we write 
\begin{align*}
   \int_0^T A_{3,1} \,\dt = &\iintTO{ \mathcal{K}_h(x-y)\lr{\pt (\Delta_x)^{-1} \dv_x [\vu]_\eta(x) - \pt (\Delta_y)^{-1} \dv_y [\vu]_\eta(y)} W }\\ 
   &+ \iintTO{ \mathcal{K}_h(x-y)\lr{\pt (\Delta_x)^{-1} \dv_x \vv_\eta(x) - \pt (\Delta_y)^{-1} \dv_y \vv_\eta(y)} W }=B_1+B_2.
\end{align*}
Note that  $W\in L^\infty(0,T;L^\infty(\T^6))$, and  $\partial_t \vu \in L^q(0,T;W^{-1,q}(\T^3))$, so for given $\eta>0$, we find easily so a certain $\alpha>0$
\begin{equation*}
    |B_1| \leq C_\eta |\ln h|^{-\alpha}, \mbox{ \  since \ } 
    \frac{1}{\|\mathcal{K}_h\|_{L^1}}\int_{\T^3} \mathcal{K}_h(z) |z| dz \sim |\ln h|^{-\alpha'}.
\end{equation*}

Next, in order to control $B_2$, first we only consider 
\begin{align*}
    \iintTO{ \mathcal{K}_h(x-y)\pt (\Delta_x)^{-1} \dv_x \vv_\eta(x) W }.
\end{align*}
Since we control the absolute value of it, we will proceed similarly for the other integral. We write the equation of $W$ as \eqh{\pt W=\dv_x \vc{A} +\dv_y \vc{B} +\mathcal{C},
}
where
\eqh{
	&\vc{A}(t,x,y)= -\vr(x) w(x)\vu(x) \sgn^\sigma_{xy} ,\\
	&\vc{B}(t,x,y)= - \vu(y) \sgn_{xy}^\sigma  \vr(x)w(x) ,\\
	&\mathcal{C}(t,x,y)= \Big(- \Lambda    \,\sgn^\sigma_{xy} - \vr(x)\dv_x\vu(x) \partial \sgn_{xy}^\sigma
	 + (\sgn^\sigma_{xy} + \vr(y) \partial \sgn^\sigma_{xy}) \dv_y \vu(y) \Big) \vr(x) w(x).}

 Using the uniform bounds w.r.t. $ \mu $, especially  $ L^\infty_{t,x} $ bound for the $w \vr$ (recall (\ref{rho-bdd})) and $L^2_{t,x}$ bound for ${\cal M}[|\Grad\vu|]$, we justify that 
\begin{align*}
	\|\vc{A}\|_{L^2_tL^6_xL^\infty_y}+\|\vc{B}\|_{L^2_tL^6_yL^\infty_x}\leq C,
\end{align*}
and
\begin{align*}
    	\Vert \mathcal{C} \Vert_{L^1_t L^1_x L^\infty_y+ L^1_t L^1_y L^\infty_x} \leq C\left( \frac{1}{\sigma},\Lambda \right)\; \Vert \vr w \Vert_{L^\infty} \left(\Vert \vr \Vert_{L^2_t L^2_x} \Vert \nabla \vu  \Vert_{L^2_t L^2_x} + \Vert \vr \Vert_{L^2_t L^2_y} \Vert \nabla \vu  \Vert_{L^2_t L^2_y} + \Vert \nabla_x \vu  \Vert_{L^2_t L^2_x}\right).
\end{align*}
Next, we introduce 
 \begin{align}
	\mathscr{W}_1(t,x)= \int_{\T^3} \mathcal{K}_h(x-y) W(t,x,y) \; \dy,\quad \mathscr{W}_2(t,y)= \int_{\T^3} \mathcal{K}_h(x-y) W(t,x,y) \; \dx . \nonumber
\end{align}
A simple observation is 
\begin{align*}
	\pt \mathscr{W}_1= & \int_{\T^3} \mathcal{K}_h(x-y) \pt W(t,x,y) \; \dy \\
	&=  \int_{\T^3} \lr{(\mathcal{K}_h(x-y) \Div_x \vc{A}(t,x,y)+\mathcal{K}_h(x-y) \Div_y \vc{B}(t,x,y) + \mathcal{K}_h(x-y) \mathcal{C}(t,x,y)) }\; \dy\\
	&= \Div_x  \lr{\int_{\T^3} \mathcal{K}_h(x-y) \vc{A}(t,x,y) \dy }+  \int_{\T^3} \mathcal{K}_h(x-y) \mathcal{C}(t,x,y) \; \dy .
\end{align*}
We have the following uniform bounds
$\mathscr{W}_1 \in  L^\infty_{t,x} \text{ with } \pt  \mathscr{W}_1 \in L^1_t W^{-1,1}_x$
such that
\begin{align*}
	&\Vert 	 \mathscr{W}_1  \Vert_{L^\infty_{t,x}} \leq \Vert \mathcal{K}_h \Vert_{L^1_x} \Vert W  \Vert_{L^\infty_{t,x,y}} ,\\
 &\Vert 	\pt \mathscr{W}_1 \Vert_{ L^2_tW^{-1,6}_x +L^1_tW^{-1,1}_x} \leq C \Vert \mathcal{K}_h \Vert_{L^1_x}.
\end{align*}
On the other hand, by definition
\begin{equation*}
\|\vve \|_{L^2(0,T;L^2)} =
\left( \int_0^T \int_{\T^3} \left( \int_{\T^3} \kappa_\eta (z) |z|
\frac{u(t,x+z)-u(t,x)}{|z|} dz 
\right)^2dx dt  \right)^{1/2} \leq C\eta \|\nabla \vu\|_{L^2(0,T;L^2)}    
\end{equation*}
and
\begin{align*}
	\pt  \vve \in  L^{q^\prime}_t W^{-1,q^\prime}_x \text{ for } q^\prime >1. 
\end{align*}
Again, recall that $\vu  $ (as a sequence $  \{\vu_\mu\}_{\mu>0}$  )  is uniformly bounded in $L^\infty_t L^2_x \cap L^2_t W^{1,2}_x $.

Here, we are in the position with the lowest possible regularity. In order to control the remaining term there is a need of support of the interpolation result:

\begin{lemma}\label{interpolation-lemma}
	Let $ \pmb{\varphi} \in L^q(\R^{d+1}) $ with $ \partial_t \pmb{\varphi} \in L^{q^\prime} (\R; W^{-1,q^\prime} (\R^d)) $ and let $ W \in L^{\bar{q}}(\R^{d+1} )  $ with $ \partial_t W \in  L^{{\bar{q}}^\prime} (\R; W^{-1,{\bar{q}}^\prime} (\R^d))  $ such that $ 1\leq q, q^\prime, \bar{q}, \bar{q}^\prime \leq \infty $ with either one of the following relations satisfied:
	\begin{align}\label{il:a1}
		\begin{split}
			&(1) \; \qquad\frac{1}{q^\prime}+\frac{1}{\bar{q}}<1<	\frac{1}{q}+\frac{1}{\bar{q}^\prime },\\
			& (2)\; \qquad\frac{1}{q^\prime}+\frac{1}{\bar{q}}>1>	\frac{1}{q}+\frac{1}{\bar{q}^\prime }.
		\end{split}
	\end{align}
	Then there exists $ 0<a<1 $ such that 
	\eq{
		\int_{\mathbb{R}} \int_{\mathbb{R}^d}  \Big((- \Delta_x)^{-1}\dv_x \pt \pmb{\varphi} \Big) W(t,x)\, \dx\; \dt 
		\leq C \Vert \pmb{\varphi}  \Vert_{L^q_{t,x}}^{a} \Vert \partial_t \pmb{\varphi} \Vert_{L^{{q}^{\prime}}_tW^{-1,{q}^\prime}_x}^{1-a}\Vert W \Vert_{L^{\bar{q}}_t}^{1-a} \Vert \partial_t W \Vert_{L^{\bar{q}^{\prime}}_t W^{-1,\bar{q}^\prime}_x}^{a}.
	}
\end{lemma}

\noindent
The complete proof of the above result is available in \cite[Lemma 15]{CMZ}. The straightforward application of Lemma \ref{interpolation-lemma} gives 
\begin{align}\nonumber
	\begin{split}
 &\iintTO{ \mathcal{K}_h(x-y)\pt (\Delta_x)^{-1} \dv_x \vv_\eta(x) W }\\
		&\leq  \left\vert \int_0^T \!\!\!\int_{\Td} (- \Delta_x)^{-1}\Div_x( \pt (\vve)) \mathscr{W}_1 \dx \; \dt \right\vert \\
		& \leq \Vert  \vve \Vert_{L_{t,x}^{q}}^{\alpha} \Vert \pt  \vve \Vert_{L^{q^\prime}_t W^{-1,q^\prime}_x}^{1-\alpha}\Vert \mathscr{W}_1 \Vert_{L^{\infty}_{t,x}}^{1-\alpha} \Vert \partial_t \mathscr{W}_1 \Vert_{L^{1}_tW^{-1,1}_x}^{\alpha}\\
		&\leq C \eta^\vt \Vert \mathcal{K}_h \Vert_{L^1_x} ,
	\end{split}
\end{align}
for some $\vt>0$, follows from the interpolation of $\vve$ and property of mollifier. In particular, we have \[ \| \vve \|_{L_{t,x}^{q}} \leq \| \vve \|_{L^2_t W^{1,2}_x}^{\alpha_1} \| \vve \|_{L^\infty_t L^2_x}^{1-\alpha_1}  \text{ with }0<\alpha_1<1.\]
Therefore, it yields 
\eq{\label{a32b}
\int_0^T A_{3,1} \, {\rm d}t 
&\leq C \|\mathcal{K}_h\|_{L^1_x} \left( |\ln h|^{-\beta} \eta^{-\alpha} +
     \eta^\vt \right).
}
Hence, collecting estimates \eqref{a-a31} and \eqref{a32b} and  $\displaystyle \widetilde{R}_{h}^{\sigma}= \frac{R_{h}^{\sigma}}{\Vert \mathcal{K}_h \Vert_{L^1_x}}$  we have 
\eq{\label{St-S0-1}
&
 \widetilde{R}_{h}^{\sigma}(T) - \widetilde{R}_{h}^{\sigma}(0)\\
&\leq C    \iintTO{  \Ov{\mathcal{K}}_{h}(x-y)\vert{\vr(y)-\vr(x)}\vert^\sigma   w(x)}\\
&\quad +C\iintTO{\Ov{\mathcal{K}}_h(x-y) \lr{D_{|x-y|}\vue(x)-D_{|x-y|}\vue(y)} \vert{\vr(x)-\vr(y)}\vert^\sigma w(x)}\\
&\quad+ C(h^{\beta}\eta^{-\alpha}+\eta^\vt)+C\sigma,
}
where $\displaystyle \Ov{\mathcal{K}}_h=\frac{{\mathcal{K}_h}}{|| \Ov{\mathcal{K}_h} ||_{L^1}}$ and $\alpha,\beta,\sigma>0$.\par
We recall the following Lemma (See, \cite[Lemma 6.3]{BJ} and \cite[Lemma 13]{CMZ}): 
\begin{lemma}\label{Lemma:Du}
Let $u \in W^{1,2}(\T^3)$ then the following holds
\begin{equation}\nonumber
\int_{\T^d} \frac{{\rm d}z}{(|z|+h)^d} \| D_{|z|}u(\cdot) - 
D_{|z|}u(\cdot + z)\|_{L^2} \leq C|\ln h|^{1/2} \|u\|_{W^{1,2}(\Td)},
\end{equation}
where $D$ is given by \eqref{formula:D}. 
\end{lemma}

Using Lemma \ref{Lemma:Du}, along with the property of $\mathcal{K}_h$ \[\|\mathcal{K}_h\|_{L^1} \sim |\ln h| \text{ as } h\approx 0^+, \]
and considering $z=x-y$, we conclude 
\eq{\nonumber
&\iintO{\Ov{\mathcal{K}}_h(x-y) \lr{D_{|x-y|}\vu(x)-D_{|x-y|}\vu(y)} |\vr(y)-\vr(x)|^\sigma w(x)}\\
&=\int_{\T^6}{\Ov{\mathcal{K}}_h(z) \lr{D_{|z|}\vu(y+z)-D_{|z|}\vu(y)} |\vr(y)-\vr(y+z)|^\sigma w(z+y)}\, {\rm d}z\,\dy\, \\
&\leq C \int_{\T^3}{\Ov{\mathcal{K}}_h(z) \|D_{|z|}\vu(\cdot+z)-D_{|z|}\vu(\cdot)\|_{L^2_x}}\, {\rm d}z\,\\
&\leq C |\log h|^{-1/2}\|\vu(t)\|_{W^{1,2}_x}.
}
Hence, \eqref{St-S0-1} reduces to 
\begin{align*}
    &\widetilde{R}_{h}^{\sigma}(T) \leq \widetilde{R}_{h}^{\sigma}(0) + \int_0^T \widetilde{R}_{h}^{\sigma}(t) \; \dt + C |\log h|^{-1/2}+ C(h^{\beta}\eta^{-\alpha}+\eta^\vt+\sigma).
\end{align*}
First we choose $\eta=h^{\frac{\beta}{2\alpha}}$, then applying the Gr\"onwall's inequality and passing the limit $\sigma \rightarrow 0$,  we get
\begin{align}
    \underset{\mu\to0}{\lim\sup} \left(  \int_0^T\!\!\!\iintO{
			{\Ov{\mathcal{K}}}_h(x-y)|\vr_\mu(t,x)-\vr_\mu(t,y)|\; w_\mu(t,x)\;} \, \dt\right) \to 0, \ \mbox{ as } h\to 0.
\end{align}
The removal of weight follows the similar line of \cite[Lemma 14]{CMZ}. 
The key idea for the removal of the weight lies in the proper use of the bound of $\vr |\log w |$ as mentioned in \eqref{eq:rhologw}. 

Let $\zeta<1$. For fixed $t$, we define $\Omega_\zeta = \{x: w^\ep(t,x) \leq \zeta\}$ and denote
$\Omega_\zeta^c$ by its complementary. 
\eq{\nonumber
&\underset{\mu\to0}{\lim\sup} \left(  \int_0^T\!\!\!\iintO{
 \Ov{\mathcal{K}}_h(x-y)|\vr_\mu(t,x)-\vr_\mu(t,y)|} \, \dt\right)\\
 &=
 \underset{\mu\to0}{\lim\sup} \left(  \int_0^T\!\!\!\int_{x\in\Omega^c_\zeta \cup y\in\Omega^c_\zeta}
 \Ov{\mathcal{ K}}_h(x-y)|\vr_\mu(t,x)-\vr_\mu(t,y)| \dx\dy\, \dt\right)\\
&\quad + \underset{\mu\to0}{\lim\sup} \left(  \int_0^T\!\!\!\int_{x\in\Omega_\zeta \cap y\in\Omega_\zeta}
 \Ov{\mathcal{K}}_h(x-y)|\vr_\mu(t,x)-\vr_\mu(t,y)| \, \dx\,\dy\, \dt\right)\\
 &=I_1+I_2.
 }
 The first term is comparably easy to compute, as we note that
\eq{
I_1&= \underset{\mu\to0}{\lim\sup} \left(  \int_0^T\!\!\!\int_{x\in\Omega^c_\zeta \cup y\in\Omega^c_\zeta}
 \Ov{\mathcal{K}}_h(x-y)|\vr_\mu(t,x)-\vr_\mu(t,y)| \dx\dy\, \dt\right)\\
 &\leq \underset{\mu\to0}{\lim\sup} \frac{1}{\zeta} \left(  \int_0^T\!\!\!\int_{x\in\Omega^c_\zeta \cup y\in\Omega^c_\zeta}
 \Ov{\mathcal{K}}_h(x-y)|\vr_\mu(t,x)-\vr_\mu(t,y)| (w_\mu(t,x){+w_\mu(t,y)})\dx\dy\, \dt\right) .
}
For the term $I_2$ we use that both $\zeta$ and $w$ are less than 1, and so $|\log w|\geq |\log\zeta|$, and so
\eq{\label{Remove-w-I2}
I_2&=  \underset{\mu\to0}{\lim\sup} \left(  \int_0^T\!\!\!\int_{x\in\Omega_\zeta \cap y\in\Omega_\zeta}
 \Ov{\mathcal{K}}_h(x-y)|\vr_\mu(t,x)-\vr_\mu(t,y)| \, \dt\right)\\
 &\leq \underset{\mu\to0}{\lim\sup} \left(  \int_0^T\!\!\!\int_{x\in\Omega_\zeta \cap y\in\Omega_\zeta}
 \Ov{\mathcal{K}}_h(x-y)|\vr_\mu(t,x)-\vr_\mu(t,y)|\frac{|\log w_\mu(t,x)|}{|\log\zeta|}  \dx\, \dy\, \dt\right)\\
 &\leq \sup_{\mu>0} \frac{2}{|\log\zeta|}\iintTO{ \mathcal{K}_h(x-y)\vr_\mu(t,x)|\log w_\mu(t,x)| }\\
 & \leq \sup_{\mu >0}  \frac{2}{|\log\zeta|} \int_0^T\!\!\!\int_{\T^3}{ \vr_\mu(t,x)|\log w_\mu(t,x)| } \dx \\
 & \leq \frac{C}{|\log\zeta|},
}
where, we applied the bound \eqref{boundlogw}.\par

Finally, using \cite[Section 8.3]{CMZ}, we relax the assumption $0\leq \vr \in L^\infty(\T^3)$ to  $\displaystyle 0\leq \vr_0 \in L^\gamma(\T^3) $. 
The key idea here is to consider the initial data for weight equation \eqref{weight-eqn}, with initial data 
\eq{\nonumber
w_{0,M}:=\begin{cases}
1 &\mbox{ for } \vr_0\leq M,\\
\frac{M}{\vr_0} &\mbox{ for } \vr_0> M,
\end{cases}
}
where, $M>>1$ be a fixed constant. Denoting, the weight as $w_M$, we note 
\begin{align}\nonumber
   \Vert \vr w_M \Vert_{L^\infty_{t,x} } \leq \Vert \vr_0 w_{0,M} \Vert_{L^\infty_x}\leq M.
\end{align} This allows us to perform the weighted convergence, as done in this section. The only difficulty we have is to perform the step where we remove the weights. The uniform bound of $ \vr \vert \log w_M \vert $ in $L^\infty_t L^1_x$ follows from 
\begin{equation*}
\underset{t>0}{{\rm ess}\sup}\intO{ \vr |\log w_M| } \leq \intO{ \vr_0 |\log w_{0,M}| } + C  \Lambda,
\end{equation*}
where $C$ is independent of $M$. Now, from the definition of $ w_{0,M} $, we obtain 
\begin{align*}
    \intO{ \vr_0 |\log w_{0,M}| }\leq \int_{\{\vr>M\}} \vr_0 |\log w_{0,M}| \dx  \leq \int_{\{\vr>M\}} \vr_0 |\log \vr_{0}| \dx \leq C_\gamma  \int \vr_0^\gamma \dx.
\end{align*}
Therefore, we conclude that 
\begin{align}
   \underset{t>0}{{\rm ess}\sup}\intO{ \vr |\log w_M| } \leq C (\vr_0,\Lambda),\nonumber
\end{align}
where $C$ is independent of $M$.

\appendix 
\section{Appendix}

\subsection{Proof of Lemma \ref{lem:evf:ep}}\label{pf:lem:evf:ep}
From the effective viscous flux equation \eqref{evf2} we note that 
\begin{align*}
	\inttauOM{G_{\ep} \; [p_{\mu}(\vr_{\ep})]_{\ep} \psi} = -\inttauOM{ \pt  (-\Delta)^{-1} \Div \vu_{\ep}  [p_{\mu}(\vr_{\ep})]_{\ep} \psi  }.
\end{align*}
Now, we introduce the space-time regularization of $\vu_\ep$ and denote it by $[\vu_\ep]_{\eta}$. We also introduce the notation \[ {\vv_{\ep}}_\eta = \vu_{\ep}- [\vu_{\ep}]_{\eta}¨\]  
and  we have $\|{\vv_{\ep}}_\eta\|_{L^{\frac{10}{3} -\lambda}_{t,x}} = o(1)$ as $\eta \to 0$, uniformly with respect to $\ep$
for any $\lambda >0$ sufficiently small.
Therefore, we write 
\begin{align*}
	& \inttauOM{  \partial_{t} {(-\Delta)}^{-1}\Div \vu_{\ep}  \;    [p_{\mu}(\vr_{\ep})]_{\ep} \; \psi } \\
	&= \inttauOM{  \partial_{t}{(-\Delta)}^{-1}\Div [\vu_\ep]_{\eta}  \;  [p_{\mu}(\vr_{\ep})]_{\ep}  \; \psi } + \inttauOM{ \partial_{t} {(-\Delta)}^{-1}\Div {\vv_\ep}_\eta  \;   [p_{\mu}(\vr_{\ep})]_{\ep}  \; \psi }=A_1+A_2.
\end{align*}
\begin{itemize}
	\item \textbf{Term $A_1$:}
	From the strong convergence of $\vu_\ep$ in $L^2_{t,x}$ and weak convergence of $ [p_{\mu}(\vr_{\ep})]_{\ep} $ in $L^{p_2}_{t,x} $ it follows for $\ep \to 0$
	\begin{align*}
		-\inttauOM{\partial_{t} {(-\Delta)}^{-1} \Div [\vu_\ep]_\eta   [p_{\mu}(\vr_{\ep})]_{\ep}  \psi  } \rightarrow -\inttauOM{\partial_{t} {(-\Delta)}^{-1} \Div [\vu]_\eta   \Ov{p_\mu(\vr)} \psi }.  
	\end{align*}  
	Now from 
	$ \pt  (-\Delta)^{-1}\Div [\vu_\ep]_\eta = -[G_\ep]_\eta$
	we have 
	\begin{align}\label{A-1-inf}
		\inttauOM{ 	\pt  (-\Delta)^{-1} \Div [\vu]_\eta \Ov{p_{\mu}(\vr)}\;\psi}& = -\inttauOM{[G]_\eta \Ov{p_{\mu}(\vr)}\; \psi} \nonumber\\
  &= -\inttauOM{G \Ov{p_{\mu}(\vr)}\; \psi} + o(1)
	\end{align} 
 for $\eta \to 0$.
	\item \textbf{Term $A_2$:}
	\begin{align}\label{a2-inf}
		&\inttauOM{ \partial_{t} {(-\Delta)}^{-1}\Div {\vv_\ep}_\eta  \;  [p_{\mu}(\vr_{\ep})]_{\ep}  \; \psi }  \nonumber \\
		&= -\inttauOM{  {(-\Delta)}^{-1}\Div {\vv_\ep}_\eta  \;   \partial_{t}( [p_{\mu}(\vr_{\ep})]_{\ep}  \; \psi) }  + \left[\intO{  {(-\Delta)}^{-1}\Div {\vv_\ep}_\eta  \;    [p_{\mu}(\vr_{\ep})]_{\ep}  \; \psi } \right]_{t=0}^{t=\tau}.
	\end{align}
	The equation for $\displaystyle ( [p_{\mu}(\vr_{\ep})]_{\ep}  \; \psi) $ reads as 
	\begin{align*}
		\partial_{t}( [p_{\mu}(\vr_{\ep})]_{\ep}  \; \psi) + \Div\left( \left[p_\mu(\vr_\ep)[\vu_\ep]_\ep \right]_{\ep}\psi \right) +[ \lr{ \vr_\ep p_{\mu}^\prime(\vr_\ep)-p_{\mu}(\vr_\ep)} [\Div \vu_\ep]_\ep ]_{\ep} \psi\\
		= -[\delta \vr_\ep^{m}p_{\mu}^{\prime}(\vr_{\ep})]_\ep\psi+ ([p(\vr_\ep)]_\ep \partial_t \psi + [p(\vr_\ep)\vu_\ep]_\ep \cdot \nabla  \psi) .
	\end{align*}
	From the uniform bounds \eqref{unif-est-ep}  we have 
	\begin{align*}
		\Vert \partial_{t} ( [p_{\mu}(\vr_{\ep})]_{\ep}  \; \psi) \Vert_{L^{s}_t W^{-1,s}_x} \leq C(\Vert \psi \Vert_{C^1_{t,x}} ).
	\end{align*}
	
We notice that the term $\displaystyle -\inttauOM{  {(-\Delta)}^{-1}\Div {\vv_\ep}_\eta  \;   \partial_{t}[p_{\mu}(\vr_\ep)]_\ep \; \psi }$ can be estimated as 
	\begin{align*}
		&\left\vert \inttauOM{  {(-\Delta)}^{-1}\Div {\vv_\ep}_\eta  \;   \partial_{t}( [p_{\mu}(\vr_{\ep})]_{\ep}  \; \psi) }  \right\vert  \\
		&= \left\vert \inttauOM{ {(-\Delta)}^{-\frac{1}{2}} \Div {\vv_\ep}_\eta    {(-\Delta)}^{-\frac{1}{2}} \partial_{t}\lr{ [p_{\mu}(\vr_{\ep})]_{\ep}  \; \psi}  }\right\vert \\
		&\leq C \Vert {\vv_\ep}_\eta \Vert_{L^{s'}_{t,x} } \Vert  {(-\Delta)}^{-\frac{1}{2}} \partial_{t} \lr{  [p_{\mu}(\vr_{\ep})]_{\ep}  \; \psi}\Vert_{L^s_{t,x} },
	\end{align*} 
	where $s>\frac{10}{7}$ and $s' < \frac{10}{3}$. Thus this term is of order $o(1)$ for $\eta \to 0$.
%
	
	The boundary term at $0$ in \eqref{a2-inf} is of the order $o(1)$ due to the choice of smooth initial data. 
	Now for the term $\displaystyle \intO{  {(-\Delta)}^{-1}\Div {\vv_\ep}_\eta(\tau)  \;   [p_{\mu}(\vr_\ep)]_\ep(\tau) \; \psi(\tau) } $ we have 
	\begin{align*}
		\left \vert\intO{ \! {(-\Delta)}^{-1}\Div {\vv_\ep}_\eta(\tau)  \;   [p_{\mu}(\vr_\ep)]_\ep(\tau) \; \psi(\tau) } \right \vert\! \leq \!C \Vert \psi \Vert_{C^1_{t,x}} \Vert  {(-\Delta)}^{-1}\!\Div {\vv_\ep}_\eta \Vert_{L^\infty_t L^{\frac{5}{3}}_x}\! \Vert[ p_\mu(\vr_\ep)]_\ep \Vert _{L^\infty_t L^\frac{5}{2}_x}\!.
	\end{align*}
	Since, $ \vu_\ep$ is bounded in $L^\infty_t L^2_x$, we have  $\displaystyle {(-\Delta)}^{-1}\Div \vu_\ep$ is bounded in $L^\infty_t W^{1,2}_x$. As $\partial_t \vu_\ep$ is bounded in $L^2_t W^{-1,2}_x$, we have that $\vu_\ep$ is bounded in $C([0,T];W^{-1,2}_x)$. Thus 
	\[\Vert  {(-\Delta)}^{-1}\Div {\vv_\ep}_\eta \Vert_{ C([0,T];L^{\frac{5}{3}}_x)} = o(1)  \quad \text{ for } \eta\to 0.\]
	Hence, combining the previous discussion we conclude that $A_2$ is of order $o(1)$ for $\eta\to 0$. 
 This finishes the proof of Lemma \ref{lem:evf:ep}.
\end{itemize}

\subsection{Existence with fixed parameter;\; Proof of Theorem \ref{exi-fix}}
\label{App.1}

{\it Proof} . It can be obtained by application of the Schauder fixed point theorem: 
\textit{``Let ${\cal T}$ be a continuous compact operator from $X$ into $X$, where $X$ is a bounded convex subset of a Banach space.  Then ${\cal T}$ has a fixed point."}\\
We assume, without loss of generality, $\pi_\mu(\vr) = p_\mu(\vr)$ as the non-monotone part of the pressure does not present any problems, except for a few technicalities. 

First we consider a bounded convex subset $X_T$ of a Banach space as 
\[\displaystyle X_T= \{ \mathbf{v} \in L^\infty(0,T; L^2(\T^3)) \cap L^2(0,T; W^{1,2}(\T^3)) \,|\, \| \mathbf{v} \|_{L^\infty_t L^2_x\cap L^2_t W^{1,2}_x} <R \} ,\]
where we will specify $R$ later.
Given $\tilde {\vu} \in X_T$, we consider the following map:
\begin{equation}\label{T}
	{\cal T}(\tilde \vu) = \vu,
\end{equation}
where $\vu $ satisfies 
\begin{align}
	&\pt \vr + \Div(\vr [\tilde \vu]_{{\ep}})+ \delta \vr^m =0,\label{bPF:1:ep-fp}\\
	&\pt \vu+ \nabla [p_\mu(\vr) ]_{\ep}= \Delta \vu \label{bPF:2:ep-fp}
\end{align}
with the initial condition $(\vu^0_\mu,\vr^0_\mu)$.
We recall the notation $[\vu]_\ep=  \kappa_\ep \ast {\tilde \vu}$ and $[p_\mu(\vr)]_\ep=  \kappa_\ep \ast  p_\mu(\vr)$.
Our strategy is to prove the desired result in the following two steps:
\begin{itemize}
	\item  First, we will prove that there exists a $T_0>0$ such that the map ${\cal T}$ defined in \eqref{T} has a fixed point.
	\item Then from the energy estimate, we recover certain bounds which will allow us to extend the solution beyond $T_0$.
\end{itemize}
Equation \eqref{bPF:1:ep-fp} can be simply solved by the method of characteristics. Moreover, one finds the following estimate (for any $p<\infty$) for the unique solution $\vr=\vr([\tilde\vu]_\ep)$ to the continuity equation:
\begin{equation*}
	\sup_{t\in (0,T)} ( \|\vr(t)\|_{L^p_x} + \|\nabla \vr(t)\|_{L^p_x}) \leq \left(\|\vr_\mu^0\|_{L^p_x}+\|\nabla \vr_\mu^0\|_{L^p_x}+\|\nabla^2 [\tilde\vu]_\ep\|_{L^1_tL^p_x}\right)
	\exp\left\{ \int_0^T \|\nabla [\tilde\vu]_\ep\|_{L^\infty_x} (t)\; \dt\right\}.
\end{equation*}
Now we turn our attention to \eqref{bPF:2:ep-fp}, which is a parabolic equation for $\vu$. The above estimate implies that $\nabla [p_\mu(\vr)]_\ep  \in L^p_tL^p_x$ for any $1<p <\infty$. Therefore, from parabolic regularity theory \cite{LSU,DHP}, for any $1<p <\infty$  we obtain unique
\begin{align*}
	\vu \in W^{1,p}(0,T; L^{p}(\T^3))\cap L^p(0,T;W^{2,p}(\T^3)).
\end{align*}
In particular, it yields the following estimates
\begin{equation*}
	\|\vu,\pt\vu,\nabla^2 \vu \|_{L^p_{t,x}} \leq C,
\end{equation*}
which will provide the compactness of the operator ${\cal T}$ in $X_T$ once we justify that ${\cal T}$ maps the ball in $L^\infty_t L^2_x\cap L^2_t W^{1,2}_x$ into itself for suitable $T>0$.

Testing \eqref{bPF:2:ep-fp} by $\vu$ we have
\begin{equation}\label{fp2}
	\frac{d}{dt} \intO{ \frac{\vert\vu\vert^2}{2}}  +\intO{|\nabla \vu|^2}  \leq \| [ p_\mu (\vr) ]_\ep \|_{L^2_x} \|\Div \vu\|_{L^2_x} \leq \| \kappa_{\ep} \|_{L^2} \| p_\mu (\vr)  \|_{L^1_x} \|\Div \vu\|_{L^2_x} .
\end{equation}
Now, multiplying the continuity equation by $\vr^{\Gamma-1}$ we get 
\begin{align}\nonumber
	\frac{d}{dt} \int_{\T^3} \vr^\Gamma \dx +  \frac{\delta \Gamma}{2}\int_{\T^3} \vr_{\ep}^{m+\Gamma-1} \dx \leq   C(\delta) \| [\Div \tilde \vu ]_{\ep} \|_{L^2_x}^2. 
\end{align}
Therefore, in particular, we have
\begin{align*}
	\| p_\mu (\vr) \|_{L^\infty_t L^1_x} \leq \|\vr_\mu^0\|^\Gamma_{L^\Gamma_x} + C(\ep,\delta)  R^2. 
\end{align*}
Plugging this into \eqref{fp2}, we derive
\begin{equation*}
	\sup_{t\leq T} \intO{ |\vu|^2(t,\cdot)} +  \iintTOM{|\nabla \vu|^2} \leq \intO{|\vu_\mu^0|^2} + 2 \|\vr_\mu^0\|^{2\Gamma}_{L^\Gamma_x} +(C(\delta,\ep)R^2 + \|\vr_\mu^0\|^\Gamma_{L^\Gamma_x})\sqrt{T} R.
\end{equation*}
We choose $\|\vr_\mu^0\|^{\Gamma}_{L^\Gamma_x} < \frac 12 \max\{R,R^2\}$, $\|\vu_\mu^0\|^2_{L^2_x}< \frac 14 R^2$ and $(C(\ep,\delta)+1)R\sqrt{T_0} <\frac 14$.

Hence we see that the ${\cal T}: X_{T_0} \rightarrow X_{T_0} $ maps a ball into itself, i.e., we have
\begin{equation*}
	\|\vu\|_{L^\infty_t L^2_x\cap L^2_t W^{1,2}_x} \leq R.
\end{equation*}
It is not difficult to verify that the operator is also continuous in the given spaces. 
Thus, using Schauder's fixed point theorem we conclude the local in time existence. To extend the local solution to the global one, we first obtain the energy estimate  
\begin{align}
	\frac{1}{2} \frac{d}{dt}\int_{\T^3}	\vert \vu \vert^2 \;\dx +   \frac{d}{dt}	\int_{\T^3} P_{{\mu}}(\vr) \;\dx + 	\int_{\T^3} \vert \nabla \vu \vert^2 \; \dx + \int_{\T^3}  \delta  \vr^{m} P_{{\mu}}^\prime(\vr) \;\dx \leq 0. \nonumber
\end{align}
This can be obtained exactly as in the proof of Lemma \ref{lem:ep:uniest}, thus we do not present details here. 
Then using the regularity of the solution and setting new initial data at $T_0$, we extend the solution beyond $T_0$.

\subsection{Additional pressure estimates}\label{bog}

We  obtain an additional estimate for pressure, in the spirit of Bogovskii estimate. We observe that we need  additional pressure estimates for both the limit passages $\delta \rightarrow 0 $ and $\mu \rightarrow 0$. 
\paragraph{Case 1: For \eqref{bPF:1:del}--\eqref{bPF:2:del}:} At first we consider the case when $\mu$ is fixed. We aim at showing that
$$
\|\vr_\delta\|_{L^{\Gamma +\alpha} ((0,T-\lambda)\times \T^3)} \leq C
$$
with $C$ independent of $\delta$ and $\alpha = \frac {13}{20} \Gamma -\frac 1{20}$ for $\Gamma$ large enough (at least $\Gamma > \frac{21}{13}$) which corresponds to the choice $m= \frac 52 \gamma +\frac 32$. Above, $\lambda >0$ arbitrarily small and $C$ blows up if $\lambda \to 0$.

In order to simplify the notation, we skip the index $\delta$ for the rest of the proof. The key idea is to consider a function $ \pmb{\varphi} $  such that 
\begin{align*}
	\dv \pmb{\varphi} = \vr^\alpha- \langle \vr^\alpha \rangle \text{ in } \T^3,
\end{align*}
where $  \langle \vr^\alpha \rangle = \frac{1}{\vert \Omega \vert} \int_\Omega \vr^\alpha \dx $ and $ \alpha $ that depends on $\Gamma$ will be chosen accordingly. However, since we deal with space-periodic boundary conditions and we want to have estimates up to the boundary, in order to avoid technicalities, we will not use directly the Bogovskii operator here and instead of it we take $\pmb{\varphi} = \nabla \psi$, where
\begin{equation} \label{Laplace}
\Delta \psi = \vr^\alpha- \langle \vr^\alpha \rangle \text{ in } \T^3
\end{equation}
(i.e., $\psi$ is also space periodic). Then we directly have
\begin{equation} \label{est_press_sob}
\|\pmb{\varphi}\|_{1,p} \leq C\|\vr^\alpha\|_p
\end{equation}
for any $1<p<\infty$. More generally, for $f \in L^{p}_x$, $\int_\Omega f \,\dx =0$, denoting $B$ the solution operator to \eqref{Laplace} with the right-hand side $f$, we have due to \eqref{est_press_sob}
\begin{equation} \label{est_press_sob1}
\|\nabla B(f)\|_{1,p} \leq C\|f\|_p
\end{equation}
and further for $f=\Div \vc{g}$, $\vc{g}$ space periodic
\begin{equation} \label{est_press_leb1}
\|\nabla B(\Div(\vc{g})\|_{p} \leq C\|\vc{g}\|_p.
\end{equation}
Multiplying the weak form of \eqref{bPF:2:del} by $\eta \nabla B(\varphi^\alpha-\langle \vr^\alpha \rangle)$, where $ \eta \in C_c^1[0,T)$, we get (we finally take $\eta =1$ for $t<T-\lambda$ for $\lambda$ small)
\begin{align*}
	-\iintOM{\eta(0)\vu_0\cdot &\nabla B\left(\vr^\alpha_0- \langle \vr^\alpha_0 \rangle \right) } =\iintTOM{ \eta\vu \cdot \partial_t \nabla B\left(\vr^\alpha- \langle \vr^\alpha \rangle \right) } \\
 &+ \iintTOM{ \eta^\prime \vu \cdot  \nabla B\left(\vr^\alpha- \langle \vr^\alpha \rangle \right) }\\
	&+ \iintTOM{\eta\pi_{\mu}(\vr) \left(\vr^\alpha- \langle \vr^\alpha \rangle \right)} 
	- \iintTOM{ \eta \nabla \vu : \nabla^2 B\left(\vr^\alpha- \langle \vr^\alpha \rangle \right) }.
\end{align*}
Using renormalized continuity equation we rewrite the above equation as 
\begin{align}\label{bog1}
	\iintTOM{\eta \pi_{\mu}(\vr) & \vr^\alpha} = -\iintOM{\eta(0)\vu_0\cdot \nabla B\left(\vr^\alpha_0- \langle \vr^\alpha_0 \rangle \right) } \\
 & +  \iintTOM{ \eta\vu \cdot \nabla B\left(\Div (\vr^\alpha\vu)\right) } \nonumber \\
 &+ \iintTOM{ \eta\vu \cdot \nabla B\left((\alpha-1)\vr^\alpha\Div \vu- \langle (\alpha-1)\vr^\alpha \Div \vu \rangle \right) } \nonumber\\
 &+ \iintTOM{ \eta^\prime \vu \cdot  \nabla B\left(\vr^\alpha- \langle \vr^\alpha \rangle \right) }
	+ \iintTOM{\eta \pi_{\mu}(\vr) \langle \vr^\alpha \rangle} \nonumber\\ 
	&+ \iintTOM{ \eta \nabla \vu : \nabla^2 B\left(\vr^\alpha- \langle \vr^\alpha \rangle \right) } + \delta \iintTOM{ \eta\vu \cdot \nabla B\left(\vr^{\alpha+m-1}- \langle \vr^{\alpha+m-1}  \rangle \right) } \nonumber\\
  &= \sum_{i=1}^{7} A_i.\nonumber
\end{align}
Below, we estimate the most restrictive terms which are $A_2$, $A_3$, $A_6$ and $A_7$. Even though the final restriction comes from the term $A_7$, we try to be precise also at other terms since this will produce the estimate in the case when $\delta \to 0$. 

\textbf{Estimate for term $A_2$:} Using \eqref{est_press_leb1} we have
\begin{align*}
	\vert A_2 \vert \leq \Vert \eta \Vert_{L^\infty_t} \Vert \vu\Vert_{L^{\frac {10}{3}}_tL^{\frac{10}{3}}_x}\Vert \vr^\alpha \vu \Vert_{L^{\frac{10}{7}}_t L^{\frac {10}{7}}_x}  \leq \Vert \eta \Vert_{L^\infty_t} \Vert \vu\Vert_{L^{\frac {10}{3}}_tL^{\frac{10}{3}}_x}^2  \Vert \vr^\alpha\Vert_{L^{\frac 52}_t L^{\frac 52}_x}. 
\end{align*}
We choose $\frac{5}{2}\alpha \leq \Gamma +\alpha$, 
i.e.,
	$\alpha \leq \frac 23\Gamma$
to obtain 
\begin{align}
	\vert A_2 \vert \leq \frac 14 	\iintTOM{\eta \pi_{\mu}(\vr)  \vr^\alpha} + C (\eta)  \Vert \vu\Vert_{L^{\frac{10}{3}}_tL^{\frac{10}{3}}_x}^{\frac {2(\Gamma+\alpha)}{\Gamma}}  .
\end{align}
\textbf{Estimate for term $A_3$:} We proceed similarly as above, however, for slightly different exponents. Recall that we know that
$$
\Vert\vu\Vert_{L^{\frac {4r}{3r-6}}_t L^r_x} \leq C, \quad 2\leq r\leq 6.
$$
We now apply \eqref{est_press_sob1} to get
\begin{align*}
	\vert A_3 \vert  \leq C \Vert \eta \Vert_{L^\infty_t} \Vert \vu \cdot \nabla  B\left(\vr^\alpha \Div \vu- \langle \vr^\alpha \Div \vu\rangle \right) \Vert_{L^1_t L^1_x} \leq \Vert \eta \Vert_{L^\infty_t} \Vert \vu  \Vert_{L^{10}_t L^{\frac{30}{13}}_x} \Vert  \Div \vu \Vert_{L^2_tL^2_x} \Vert\vr^\alpha\Vert_{L^{\frac 52}_t L^{\frac 52}_x}.
\end{align*}
We bound this term similarly as above assuming that $\alpha \leq \frac 23\Gamma$ holds.

%
\textbf{Estimate for term $A_6$:}  
We get
\begin{align*}
	\vert A_6 \vert \leq \Vert \eta \Vert_{L^\infty_t} \Vert \nabla \vu : \nabla^2 B\left(\vr^\alpha- \langle \vr^\alpha \rangle \right) \Vert_{L^1_t L^1_x} \leq \Vert \eta \Vert_{L^\infty_t} \Vert \nabla \vu  \Vert_{L^{2}_t L^{2}_x} \Vert \vr^\alpha  \Vert_{L^\infty_t L^2_x}.
\end{align*}

With additional hypothesis 
	$2\alpha \leq \Gamma +\alpha,$
i.e., $\alpha \leq \Gamma$ 
we derive 
\begin{align*}
	\vert A_6 \vert  \leq \frac 14 	\iintTOM{\eta \pi_{\mu}(\vr)  \vr^\alpha} + C(\eta)    \Vert  \nabla \vu \Vert_{L^{2}_t L^{2}_x}^{\frac{\Gamma+\alpha}{\Gamma}}.
\end{align*}
\textbf{Estimate for term $A_7$:} We proceed as above. 
\begin{align*}
	\vert A_7 \vert \leq \delta \Vert \eta \Vert_{L^\infty_t} \Vert \vu \cdot  \nabla B\left(\vr^{\alpha+m-1}-\langle \vr^{\alpha+m-1}\rangle \right) \Vert_{L^1_t L^1_x} \leq \delta \Vert \eta \Vert_{L^\infty_t} \Vert \vu  \Vert_{L^{10}_t L^{\frac{30}{13}}_x} \Vert \vr^{\alpha+m-1}  \Vert_{L^{\frac {10}{9}}_t L^{\frac{10}{9}}_x}.
\end{align*}
We observe that 
$\frac{10}{9}(m-1+\alpha)\leq m-1+ \Gamma$,
provided
$m-1 \leq 9 \Gamma -10\alpha.$
Thus we can take $\alpha = \tfrac {13}{20} \Gamma -\tfrac 1{20}< \frac 23 \Gamma$ which corresponds to our choice of $m= \frac 52\Gamma +\frac 32$. To have $\alpha >1$, we  need $\Gamma > \frac{21}{13}$.

Hence, we deduce that 
\begin{align*}
	\|\pi_\mu(\vr_\delta) \vr^\alpha\|_{L^1((0,T-\lambda)\times \T^3)} \approx \|\vr_\delta^{\Gamma+\alpha}\|_{L^1((0,T-\lambda)\times \T^3)} \leq C.
\end{align*}
This yields that there exists $p>1$ such that we have 
$ \|\pi_\mu(\vr_\delta)\|_{L^p((0,T-\lambda)\times \T^3)} \leq C.$

\paragraph{Case 2: For  \eqref{bPF:1:ep}--\eqref{bPF:2:ep}:}
We now consider $\mu \to 0$. We intend to show that
$$
\|\vr_\mu\|_{L^{\gamma +\alpha} ((0,T-\lambda)\times \T^3)} \leq C
$$
with $C$ independent of $\delta$ and $\alpha \leq \frac 23 \gamma$.
As above, we skip below the index $\mu$.  Proceeding with very similar calculation as above we obtain \eqref{bog1} without the term $\mathcal{A}_7$. Moreover, 
here want to derive a $\mu$ independent estimate, so we write the terms as 
\begin{align}\label{bog2}
	\begin{split}
		\iintTOM{\eta \pi(\vr) &\vr^\alpha} =  -\iintOM{\eta(0)\vu_0 \cdot \nabla  B\left(\vr^\alpha_0- \langle \vr^\alpha_0 \rangle \right) }
		+\iintTOM{ \eta\vu   \cdot\nabla B\left(\Div \vr^\alpha \vu \right) }\\
  &+\iintTOM{ \eta\vu   \cdot \nabla B\left((\alpha-1) \vr^\alpha \Div \vu -\langle (\alpha-1)\vr^\alpha \Div\vu\rangle  \right) } \\
		&+ \iintTOM{ \eta^\prime \vu \cdot  \nabla B\left(\vr^\alpha- \langle \vr^\alpha \rangle \right) } +  \iintTOM{\eta\pi_\mu(\vr) \langle \vr^\alpha \rangle } \\
		&+ \iintTOM{ \eta\nabla \vu : \nabla^2 B\left(\vr^\alpha- \langle \vr^\alpha \rangle \right)}	- \mu \iintTOM{\eta\vr^{\Gamma +\alpha} } \\
		&= \sum_{i=1}^{7} A_i.
	\end{split}
\end{align}
We can ignore the term $\displaystyle - \mu \iintTOM{\eta\vr^{\Gamma +\alpha} } $ since it has good sign.
We proceed similarly as above and recalling the asymptotics of the pressure we get exactly the same results as in the previous case, ignoring the limitation from $A_7$ and replacing $\Gamma$ by $\gamma$. For
\begin{align}\label{BE-mu1}
	\alpha \leq \frac{2}{3}\gamma
\end{align}
we have 
	$\|\pi(\vr_\mu) \vr^\alpha\|_{L^1((0,T-\lambda)\times \T^3)} \approx \|\vr_\mu^{\Gamma+\alpha}\|_{L^1((0,T-\lambda)\times \T^3)} \leq C.$

\subsection{Boundedness of oscillation defect measure. Proof of Proposition \ref{prop:OscD}}
\label{a5}

\textit{Proof:} First we observe that
\begin{align}\label{def-m-f1}
	\vert T_k(\vr_\mu) -T_k(\vr) \vert^{\gamma+1} \leq (\vr_\mu^\gamma - \vr^\gamma) (T_k(\vr_\mu) -T_k(\vr)).
\end{align}
The above statement follows from the following observations:
\begin{align*}
	&\vert T_k(t) -T_k(s) \vert \leq \vert t-s \vert, \text{ for } t,s>0,\\
	&(t-s)^\gamma \leq (t^\gamma -s^\gamma),\;\text{ for } t>s>0,\\
	&\vert  T_k(t) -T_k(s)\vert^{\gamma+1} \leq (t^\gamma -s^\gamma)( T_k(t) -T_k(s)) , \; t,s>0.
\end{align*}
The fact 
	$\vr \mapsto \vr^\gamma \text{ is convex and } \vr \mapsto T_k(\vr) \text{ is concave} $
yields
\begin{align*}
	\vr^\gamma \leq \overline{\vr^\gamma}\text{ and } T_k(\vr)\geq \overline{T_k(\vr)}.
\end{align*}
Hence, we have 
\begin{align*}
	\limsup_{\mu\to 0} 	\iintTOM{\vert T_k(\vr_\mu) -T_k(\vr) \vert^{\gamma+1} }&\leq \iintTOM{(\ \overline{\vr^\gamma T_k(\vr)} - \overline{\vr^\gamma} \overline{T_k(\vr)} )}\\
	&\quad + \iintTOM{\underbrace{\left( ( \overline{\vr^\gamma}- \vr^\gamma) (\overline{T_k(\vr)}-T_k(\vr))\right)}_{\leq 0}}.
\end{align*}

We recall the effective viscous flux identity with the truncation function $T_k$:
\begin{align*}
	\overline{\pi(\vr) T_k(\vr) }- \overline{\pi(\vr)}\;\overline{T_k(\vr)} =\overline{T_k(\vr)\dv \vu} - \overline{T_k(\vr)} \dv \vu .
\end{align*}
We now distinguish two cases.
\subsection*{Monotone case $(\pi(\vr) =a \vr^\gamma + \tilde p(\vr)$):}
Therefore, using structure of $\pi $ ($\tilde p$ is monotone), we have 
\begin{align*}
&a\overline{ \vr^\gamma T_k(\vr)} - a \overline{\vr^\gamma} \; \overline{T_k(\vr)} \leq a\overline{ \vr^\gamma T_k(\vr)} + \overline{\tilde p(\vr)T_k(\vr)} -a \overline{\vr^\gamma}\; \overline{T_k(\vr)} - \overline{\tilde p(\vr)}\; \overline{T_k(\vr)}\\
& =	\overline{\pi(\vr) T_k(\vr) }- \overline{\pi(\vr)}\;\overline{T_k(\vr)} =\overline{T_k(\vr)\dv \vu} - \overline{T_k(\vr)} \dv \vu =\lim\limits_{\mu \to 0} \lr{(T_k(\vr_\mu) -  \overline{T_k(\vr)} ) \dv \vu_\mu}. 
\end{align*}
Above, we used \cite[Theorem 11.26]{FN} yielding  $\overline{\tilde p(\vr)T_k(\vr)} \geq \overline{\tilde p(\vr)}\; \overline{T_k(\vr)} $.
The computations above result into
\begin{align*}
	\limsup_{\mu\to 0} 	\iintTOM{\vert T_k(\vr_\mu) -T_k(\vr) \vert^{\gamma+1} } & \leq C\limsup_{\mu\to 0} 	\iintTOM{ \left\vert \lr{(T_k(\vr_\mu) -  \overline{T_k(\vr)} ) \dv \vu_\mu} \right\vert }  \\
	& \leq C\left \Vert  \dv \vu_\mu \right\Vert_{L^2_{t,x}}\limsup_{\mu\to 0}  \left\Vert T_k(\vr_\mu) -  \overline{T_k(\vr)}  \right\Vert_{L^2_{t,x}} .
\end{align*}
From the observation
\begin{align*}
	\limsup_{\mu\to 0} \Vert T_k(\vr_\mu) -  \overline{T_k(\vr)} \Vert_{L^2((0,T)\times \Omega)}& \leq 	\limsup_{\mu\to 0} \Vert T_k(\vr_\mu) -  {T_k(\vr)} \Vert_{L^2((0,T)\times \Omega)} \\
	&\quad \quad +	\limsup_{\mu\to 0} \Vert T_k(\vr) -  \overline{T_k(\vr)} \Vert_{L^2((0,T)\times \Omega)} \\ &\leq 2 	\limsup_{\mu\to 0} \Vert T_k(\vr_\mu) -  {T_k(\vr)} \Vert_{L^2((0,T)\times \Omega)} 
\end{align*}
together with $ \gamma+1 >2  $, we conclude 
\begin{align*}
	\limsup_{\mu\to 0} 	\iintTOM{\vert T_k(\vr_\mu) -T_k(\vr) \vert^{\gamma+1} } < \infty .
\end{align*}
\subsection*{Non monotone case with $(\pi(\vr) = p(\vr) + q(\vr) )$:}
The pressure $\pi(\vr) = p(\vr)+ q(\vr) $ can be rewritten as 
\begin{align*}
	\pi(\vr) = a\vr^\gamma + \tilde{p}(\vr) + \tilde{q}(\vr)
\end{align*}
where $ \tilde{p}(\vr) $ is a monotone function and $\tilde{q}$ is uniformly bounded (as well as compactly supported). 
Then, following similar computations as before, we end up with 
\begin{align*}
	&a \limsup_{\mu\to 0} 	\iintTOM{\vert T_k(\vr_\mu) -T_k(\vr) \vert^{\gamma+1} } \\
	&\leq (\left\Vert  \dv \vu_\mu \right\Vert_{L^2_{t,x}}+ C(\tilde q))\limsup_{\mu\to 0}  \left\Vert T_k(\vr_\mu) -  \overline{T_k(\vr)}  \right\Vert_{L^2_{t,x}}  
\end{align*}
and we conclude as above.

\subsection{Bounded oscillation defect measure implies renormalized continuity equation. \; Proof of Proposition \ref{prop:REC}}
\label{a6}

\textit{Proof:} As $T_k(\vr_\mu) \to \overline{T_k(\vr)}$ in $C_{weak}([0,T]; L^\beta(\T^3))$ for any $1\leq \beta<\infty$ and $T_k(\vr_\mu) \vu_\mu \to \overline{T_k(\vr)}\vu$ in $L^2((0,T)\times \T^3)$, passing to the limit $\mu \to 0$ in the renormalized continuity equation we end up with
$$
(\overline{T_k(\vr)})_t + \dv (\overline{T_k(\vr)}\vu) + \overline{(T_k'(\vr) \vr-T_k(\vr))\dv \vu} =0
$$
in $\mathcal D'((0,T)\times \T^3$. We now renormalize this equation (the Friedrich commutator lemma can be clearly used in this situation) to obtain
$$
(b(\overline{T_k(\vr)}))_t + \dv (b(\overline{T_k(\vr)})\vu) + (b'(\overline{T_k(\vr)}) \overline{T_k(\vr)}  - b'(\overline{T_k(\vr)}))\dv\vu   = - b'(\overline{T_k(\vr)}) \overline{(T_k'(\vr) \vr-T_k(\vr))\dv \vu} 
$$
in $\mathcal D'((0,T)\times \T^3$, where $b\in C^1[0,\infty)$ such that $b'(z) =0$ for $z\geq M$ for some $M>0$. Since
$\lim_{\mu \to 0} \|\vr_\mu-T_k(\vr_\mu)\|_{L^1(0,T)\times \T^3)} = \|\vr-\overline{T_k(\vr)}\|_{L^1((0,T)\times \T^3)}$ and the right-hand side tends to zero when $k\to \infty$, it is enough to verify that
$$
\lim_{k\to \infty} b'(\overline{T_k(\vr)}) \overline{(T_k'(\vr) \vr-T_k(\vr))\dv \vu} =0
$$
in $L^1((0,T)\times\T^3)$. For $Q_{k,M}:= \{(t,x)\in (0,T)\times \T^3; \overline{T_k(\vr)} \leq M\}$, we have
\begin{align*}
\|b'(\overline{T_k(\vr)}) \overline{(T_k'(\vr) \vr-T_k(\vr))\dv \vu}\|_{L^1(Q_{k,M})} \\
\ \ \ \ \  \leq C\sup_{\mu} \|\dv\vu_\mu\|_{L^2((0,T)\times\T^3 )} \liminf_{\mu\to 0} \|T_k(\vr_\mu)-T_k'(\vr_\mu)\vr_\mu\|_{L^2(Q_{k,M})}.   
\end{align*}
But
$$
\|T_k(\vr_\mu)-T_k'(\vr_\mu)\vr_\mu\|_{L^2(Q_{k,M})} \leq \|T_k(\vr_\mu)-T_k'(\vr_\mu)\vr_\mu\|^\alpha_{L^1(Q_{k,M})} \|T_k(\vr_\mu)-T_k'(\vr_\mu)\vr_\mu\|^{1-\alpha}_{L^r(Q_{k,M})}
$$
for some $0<\alpha<1$. 
Since $0\leq T_k'(\vr_\mu)\vr_\mu \leq T_k(\vr_\mu)$, we get (we use the weak lower semicontinuity of the norm)
\begin{align*}
&\|T_k(\vr_\mu)-T_k'(\vr_\mu)\vr_\mu\|_{L^{\gamma+1}(Q_{k,M})}  \\
& \leq  \|T_k(\vr_\mu)-T_k(\vr)\|_{L^{\gamma+1}(Q_{k,M})} + \|T_k(\vr)-\overline{T_k(\vr)}\|_{L^{\gamma+1 }(Q_{k,M})} + \|\overline{T_k(\vr)}\|_{L^{\gamma+1}(Q_{k,M})} \\
& \leq \|T_k(\vr_\mu)-T_k(\vr)\|_{L^{\gamma+1}(Q_{k,M})} + {\rm osc}_{\gamma +1} (\vr_\mu-\vr) + M(T|\T^3|)^{\frac 1{\gamma+1}}\\
& \leq 2 {\rm osc}_{\gamma +1} (\vr_\mu-\vr) + M(T|\T^3|)^{\frac 1{\gamma+1}}.
\end{align*}
Thus also
$$
\lim_{k\to \infty} \sup_{\mu} \|T_k(\vr_\mu)-T_k'(\vr_\mu)\vr_\mu\|_{L^2(Q_{k,M})}=0.
$$
\medbreak
{\bf Acknowledgment:} The work of N.C. was supported by the ``Excellence Initiative Research University (IDUB)'' program at the University of Warsaw. The second author (P.B.M.) has been partly supported by the Narodowe
Centrum Nauki (NCN) Grant No. 2022/45/B/ST1/03432 (OPUS). The work of the third author (M.P.) was partly supported by the Czech Science Foundation, Grant No: 22-01591S.


\end{document}